\documentclass{article}
\usepackage{times}
\font \ssym = msam7	

\font \sym = msam10

\input epsf

\newtheorem{Thm}{Theorem}
\newtheorem{Def}{Definition}
\newtheorem{Lem}{Lemma}

\newtheorem{Cor}{Corollary}

\newcommand{\vs}{\vspace{5 mm}}

\def\dia{{\sym\char'006} }
\def\diaa{{\sym\char'006}$\!$}
\def\diaas{{\sym\char'006}\nobreak\hskip -0.3 pt s }
\def\diab{{\sym\char'007} }
\def\diabb{{\sym\char'007}$\!$}
\def\diabs{{\sym\char'007}\hskip -0.3 pt s }
\def\diabspecial{{\sym\char'007}\hskip -0.3 pt s\break}
\def\whbox{{\ssym\char'003} }
\def\blbox{{\ssym\char'004} }
\def\whboxx{{\ssym\char'003}\hskip -0.2 pt}
\def\blboxx{{\ssym\char'004}\hskip -0.1 pt}
\def\blboxs{{\ssym\char'004}\hskip -0.3 pt s\break}

\def\wvw{\raise2.3pt\hbox{$\mathchoice\sqr89\sqr98\sqr{2.1}4\sqr{1.5}4$}}

\begin{document}

\title{Isonemal Prefabrics with No Axes of Symmetry\thanks{Work on this material has been done at home and at Wolfson College, Oxford.
Richard Roth helped with the understanding of his papers and with this presentation. 
The workshop audience at the Banff International Research Station in April, 2004, gave valuable comments on the first diagrams.
Will Gibson made it possible to draw these diagrams with surprising ease from exclusively keyboard input using Maple and Xfig. 
An anonymous referee for this journal made many helpful suggestions.
To them all I make grateful acknowledgement.}}
\author{R.~S.~D.~Thomas
}
\maketitle

\noindent St John's College and Department of Mathematics, University of Manitoba,

\noindent Winnipeg, Manitoba  R3T 2N2  Canada.
thomas@cc.umanitoba.ca

\noindent Other contact information: phone 204 488 1914, fax 204 474 7611.

\begin{abstract} \noindent This paper refines Richard Roth's taxonomy of isonemal weaving designs through the final types 33--39 in order to complete the solution of three problems for those designs: which designs exist in various sizes, which prefabrics can be doubled and remain isonemal, and which can be halved and remain isonemal. These types have no symmetry axes but have quarter-turn symmetries. Jean Pedersen's problem of woven cubes is also discussed.
\end{abstract}

\noindent AMS classification numbers: 52C20, 05B45.

\smallskip
\noindent Keywords: fabric, isonemal, prefabric, weaving.

\vs
\noindent  {\bf 1. Introduction} 

\noindent The purpose of this paper is to refine, for prefabrics with quarter-turn symmetry but no axes of symmetry, the taxonomy of prefabrics proposed by Richard L.~Roth [7] and to use this final refinement to complete answers to a few questions about fabrics (Sections 5--8).
Together with [10] and [11], this completes a survey of all the prefabrics in Roth's 39 infinite families, exhausting all prefabrics except for a short list of interesting exceptions of small order (one-dimensional period).
As Roth observes beginning his subsequent paper [8] on perfect colourings, `[r]ecent mathematical work on the theory of woven fabrics' begins with Gr\"unbaum and Shephard's [2], which remains the fundamental reference.
Roth's paper [7], however, contains the major advance from the fundamental work of Gr\"unbaum and Shephard [2, 3, 4, 5].
In it he determines the various (layer, similar to crystallographic) symmetry-group types that periodic isonemal fabrics---actually prefabrics too---can have and, in the sequel, which of them can be perfectly coloured by striping warp and weft.
We are not concerned with striping, but the other terms are defined in [10], to which reference needs to be made.

Since our prefabric layers meet at right angles and the symmetry groups with which we shall be concerned here are all rotational, the {\it lattice} units, that is, units whose vertices are all images of a single point under the action of the translation subgroup, are well-defined, but not uniquely defined, as squares.
The standard reference for symmetry terms is [9], except that {\it order} is used here for 1-dimensional period to distinguish it from what will be called the {\it period}, the 2-dimensional period, which often differs.
The notion of symmetry group allows the definition of the term isonemal; a prefabric is said to be {\em isonemal} if its symmetry group is transitive on the strands.
\begin{figure}
\begin{tabbing}

\noindent
\epsffile{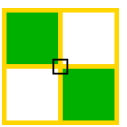}\hskip 10 pt\=
\epsffile{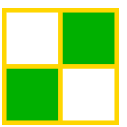}\hskip 10 pt\=
\epsffile{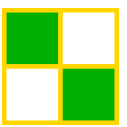}\hskip 10 pt\=
\epsffile{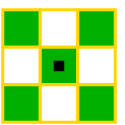}\hskip 10 pt\=
\epsffile{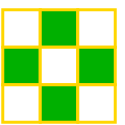}\\

\noindent (a)\>(b)\>(c)\>(d)\>(e)\\
\end{tabbing}\vskip - 18 pt

\noindent Figure 1. a. Fragment of plain weave. \hskip 10 pt b. Reverse of the fragment of Figure 1a as seen in mirror. \hskip 10 pt c. Reverse of the fragment of Figure 1a as seen from behind. \hskip 10 pt d. Fragment of plain weave. \hskip 10 pt e. Reflection of the reverse of the fragment of Figure 1d.
\end{figure}

\smallskip
If a finite region of the plane $E$ of a prefabric is viewed from behind by setting up a mirror beyond it and looking at $E$'s reflection, then what is seen in each cell (bounded in pale grey in figures) is the strand perpendicular in direction and opposite in colour in corresponding positions, as in Figure 1b compared with the front view of Figure 1a (ignoring the mark).
Figure 1b is in contrast to what one would see from the other side of $E$, which appears in Figure 1c, where the correspondence of cells of $E$ between Figure 1c and Figure 1a is obscured by the left-right reversal that causes mirror-image confusion in the real world.
For the sake of the correspondence, the reverse of a fabric will be represented as in Figure 1b; the mirror here being a simplifying device.
As long as the strands are coloured normally (warps dark, wefts pale), the reverse (so viewed) is just the colour-complement (opposite) of the obverse%
\footnote{I am using the words `obverse' and `reverse' rather than `front' and `back' because of the arbitrariness of which is which.}
 and so is of no visual interest.

Because one of the isometries that is used in weaving symmetries is reflection in the plane $E$, the reversal of which strand is uppermost at every non-boundary point, it is good to have an easy way to represent such reflection $\tau$.
As Figures 1a and 1b make clear, reversal of dark and pale represents the action $\tau$ adequately.
And that makes clear why $\tau$ can be {\it in} a symmetry but cannot itself {\it be} a symmetry.
If Figure 1a is rotated a quarter-turn about its centre, indicated by the tiny black square at that point, the pattern changes to that of Figure 1b, but, because warps become wefts under a quarter-turn and vice versa, Figure 1b is not the design of the rotated fabric. 
To be the design of the rotated fabric, the colours must be reversed, which reproduces Figure 1a, correctly indicating that that quarter-turn (no $\tau$) is a symmetry of that fabric fragment and of the whole design.
A quarter-turn without $\tau$ is represented and located by the hollow square \whboxx .
If in contrast Figure 1d is rotated a quarter-turn about its centre, the pattern remains Figure 1d, but, because warps and wefts are interchanged under a quarter-turn, Figure 1d is not the design of the rotated fabric.
The correct representation of the rotated fabric is Figure 1e, correctly indicating that a quarter-turn with the indicated centre is not a symmetry of the fabric.
If, however, $\tau$ is added to the operation, then the colours are reversed, Figure 1d is restored, and the combination of quarter-turn and $\tau$ is seen to be a symmetry of the fabric because it preserves its design.
That combination at that location is represented by the filled square \blboxx .
These notations are consistent with those in [11], where centres of half-turns are represented as \dia and of half-turns with $\tau$ as \diabb .
In later figures, black solid and dashed lines denote boundaries of lattice units or other regions of interest.

Most of the crystallographic group types are not relevant to this paper, and those that are, $p4$ and its subgroup $p2$, will be well exemplified by the fabrics to which they apply.
Symmetries involving groups of type $p4$ will be called {\it rotational} despite the existence of half-turns {\it but not quarter-turns} in symmetry groups of Roth types 11--32 [11].
The rotational symmetry groups are Roth types 33 to 39; types 1 to 32 have no quarter-turn symmetry but instead have various reflection and glide-reflection symmetries.
The small number of prefabrics of Roth types beyond 39 have both quarter-turn and reflective symmetries; Roth duly terms them {\it exceptional}.
A $p4$ lattice unit has to be square, and centres of quarter turns appear at its corners and centre.
At the mid-points of its sides there are centres of half-turns, located by \dia or \diabb .
The relation between quarter-turns and half-turns can be symbolized with slight abuse of notation as
$\hbox{\whboxx\hskip 0.7 pt}^2 = \hbox{\blboxx\hskip 0.7 pt}^2 = \hbox{\diaa} .$
That is a partial description of the groups of crystallographic type $p4$; more will follow.
Roth distinguishes seven types of these symmetry groups, one with three unnamed subtypes.
This list will be refined here from nine to eleven subtypes.
The full symmetry group of a prefabric is denoted $G$ by Roth, and within it its {\it side-preserving subgroup} is denoted $H$.

Associated with $G$ is a planar group, its planar projection denoted $G_1$ by Roth, which is algebraically isomorphic 
to $G$ and consists of all $G$'s elements $g$ whether paired with 
$\tau$ or with $e$, the identity of the group of reflections in $E$. 
$G$'s {\it side-preserving subgroup} $H$, consisting of the 
elements $(g, e)$ and omitting $(g, \tau)$ elements is also 
isomorphic to its planar projection 
$H_1$ consisting of those elements $g$ that are paired with $e$ in $H$.
$H_1$ will be referred to as the side-preserving subgroup of $G_1$.
\begin{figure}
\begin{tabbing}

\noindent\epsffile{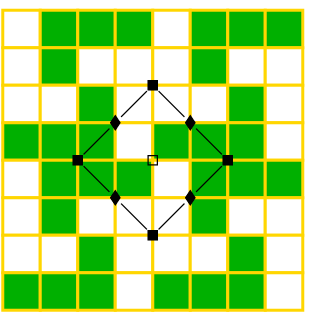}\=\epsffile{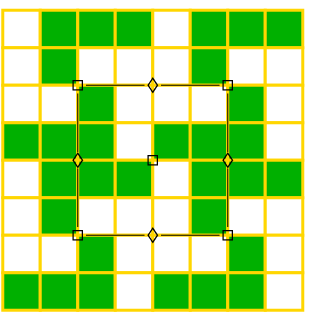}\\

\noindent (a)\>(b)\\
\end{tabbing}\vskip - 18 pt

\noindent Figure 2. The simple houndstooth, the design of the prefabric 4-1-2*, which falls apart. a. With $G_1$ displayed. \hskip 10 pt b. With $H_1$ displayed, \blbox being replaced by \dia because in $H_1$ \blbox occurs only squared.
\end{figure}

The type of a symmetry group is presented by Roth in a Coxeter notation, $G_1$/$H_1$, where the types of the groups are substituted, giving $p4/p2$, or $p4/p4$ if $H_1\neq G_1$ is of the same crystallographic type as $G_1$, or $p4/-$ if $H_1=G_1$.
Figure 2a illustrates $G_1$ for the prefabric numbered 4-1-2* in the catalogue of isonemal prefabrics that fall apart [5], the simplest houndstooth, and Figure 2b illustrates its $H_1$.
Note that the lattice unit for $H_1$ is larger than that for $G_1$, because $H_1$ is a subgroup of $G_1$.
A lattice unit for $G_1$ of a prefabric will often be referred to as `the' lattice unit of the prefabric despite its non-uniqueness and despite $H_1$'s also having a lattice unit.
Non-uniqueness will be emphasized in due course.
4-1-2* is the only isonemal prefabric with rotational symmetry in the catalogue [5] of those up to order 16 that fall apart.

The plan of the paper is first to determine what lattice units might be possible for $p4$-type symmetry groups of isonemal prefabrics. 
In \S 3 the actual lattice units of groups of the various Roth types are found and square satins in particular are discussed.
In \S 4 it is shown that isonemal prefabrics cannot have lattice units beyond those found in \S 3.
In \S\S 5--7, the announced problems are solved.
Then in \S 8 which species of prefabric can be used to weave cubes are determined, subject to some constraints.

\vs
\noindent {\bf 2. Levels of Lattice Units}

\noindent The first task is to determine which lattice units are feasible for $p4$ groups.
As Figure 2 indicates, the lattice units for prefabric 4-1-2* are conformed to the directions of the strands.
The figure cannot indicate it, but this conformity is not feasible for prefabrics of order greater than 4.
It was pointed out in [3] that prefabrics with axes of reflection or glide-reflection parallel to the strand boundaries have to have order 1, 2, or 4.
But a similar result is obviously true for isonemal prefabrics with rotational symmetry groups too, for if the lattice unit in Figure 2a is enlarged there will be no half-turns to relate many adjacent strands.%
\footnote{Since Roth has shown that every non-exceptional isonemal fabric (or prefabric) with quarter-turn symmetry is of genus III or V (in the language of [3]), although it may {\it also} be of genus I, II, or IV, every strand is transformed into its immediate neighbours by a half-turn.}
Figure 2b represents the minimal enlargement and illustrates the point; nothing relates a strand through the lattice unit to one {\em or} the other adjacent strand (whether vertical or horizontal).
This perfectly legitimate transformation group cannot be the symmetry group of an isonemal prefabric; this is a situation that will recur.
What is needed for the lattice unit to be feasible is that it be oblique to the strands, that is, that its boundaries not be parallel to or at $45^\circ$ to the strand directions.
Another feature of conformed lattice units is that the rest of the units in the plane add nothing to the whole group's capacity to interchange a particular pair of neighbouring strands.
On the other hand, for oblique lattice units centres of half-turns in a large number of different units can contribute to the interchange of strands passing through a single unit.

In order to accommodate a cell-corner quarter-turn at one corner of the lattice unit and a half-turn at the mid-point of a lattice unit's side on opposite edges of a strand, the lattice-unit side must rise two cells along its length (Figure 3a).
On the other hand, to accomodate a quarter-turn centre in the middle of a strand (and therefore of a cell) and a half-turn centre at the mid-point of the lattice unit's side on an edge of the strand, the lattice-unit side must rise one cell along its length (Figure 3b).
Other orientations are possible, and we shall get to them, but these are the simplest and illustrate that the obliquity of the lattice unit distributes half-turn centres.
It is done still more simply by placing the corners of the lattice unit on opposite sides of the same strand (Figure 3c).
This has the consequence that the work of transformation of strand to strand is done by the corners, and the half-turns at mid-side fall inside cells, giving the warp and weft the half-turn symmetry that is also present in Figure 3b in both warp and weft at \blbox as \blboxx\hskip 0.8 pt${}^2$ and in warp only at \diaa .
\begin{figure}
\begin{tabbing}

\noindent\epsffile{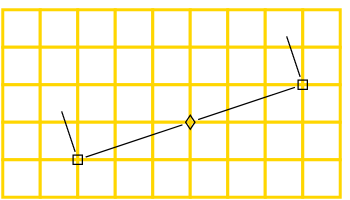}\hskip 10 pt\=\epsffile{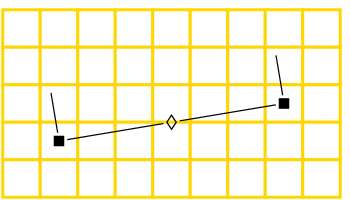}\hskip 10 pt\=\epsffile{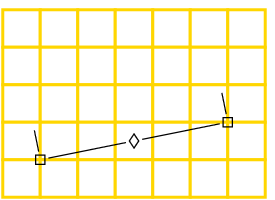}\\

\noindent (a)\>(b)\>(c)\\
\end{tabbing}
\vskip -18 pt
\noindent Figure 3. One side of a $p4$ lattice unit. a. With mid-side \dia at a cell corner. \hskip 10 pt b. With mid-side \dia on cell side. \hskip 10 pt c. With mid-side \dia within cell.
\end{figure}

We shall see both how to make lattice units as large as we please by appropriate obliquity and how to make them have area an even number of cells by doubling area, and so we begin with the smallest possible oblique lattice unit having area an odd number of cells.
By the theorem of Pythagoras, the unit's side is the hypotenuse of a right triangle with other sides of opposite parity. 
Examples (Figures 6--9) show that the centre of the unit can be in either of the two logically possible positions, cell centre and cell corner.
For definiteness, consider first the centre of a cell as in Figure 4a.
Let the origin of a co-ordinate system $(0,0)$ be placed at the centre, the axes be set in the directions of the strands, and the cell side be used as the unit along the axes. 
Keeping corners of lattice units off the lines $x=0$, $y=0$, and $y=\pm x$ ensures obliquity.
Because we shall have to refer to these lines, I call them the {\it forbidden lines}.
Using for lattice-unit corners points of cell-centre/cell-corner type opposite to the unit centre ensures oddness of unit area.
Then corners of lattice units can fall at any points off the forbidden lines and of the cell-centre/cell-corner sort opposite to $(0,0)$, i.e., $(k + {1\over2}, \ell + {1\over2})$ ($k$ and $\ell$ integers).
For definiteness and with no loss of generality, consider $k > -\ell > 0$.
For these points in the fourth quadrant, the lattice unit is the square on a hypotenuse running from
$(k + {1\over2}, \ell + {1\over2})$ to its image under a quarter-turn about $(0,0)$, $(-(\ell + {1\over2}), k + {1\over2})$, the other sides of the triangle running from $(k + {1\over2}, \ell + {1\over2})$ to 
$(k + {1\over2}, k + {1\over2})$ and from there to
$(-(\ell + {1\over2}), k + {1\over2})$ with lengths 
$$k-\ell \equiv M_1  \ \hbox{and}\  k + \ell +1 \equiv N_1.\eqno(1)$$
Any choice of parity of $k$ and $\ell$ gives differing parities to these lengths and so makes the areas of the lattice units, $M_1^2 + N_1^2$, be odd numbers starting with 5, 13, 17.
We shall regard these lattice units as basic.
\begin{figure}
\begin{tabbing}

\noindent\epsffile{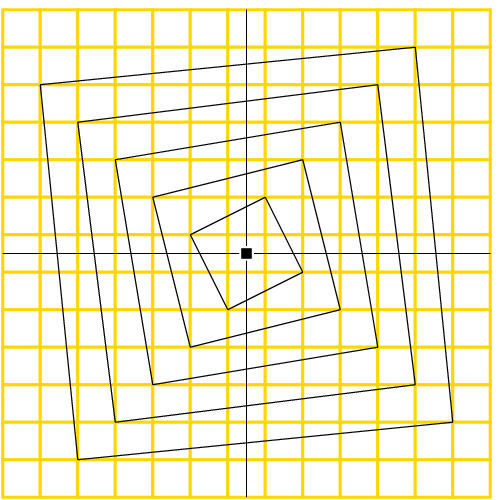}\hskip 5 pt\=\epsffile{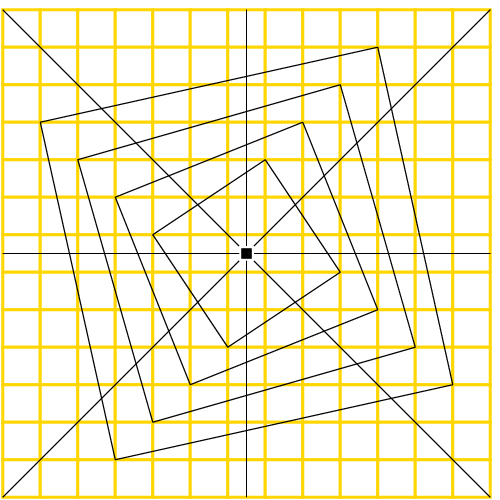}\\

\noindent (a)\>(b)\\
\end{tabbing}
\vskip -18 pt
\noindent Figure 4. Two arbitrary growing sequences of level-1 lattice units with centre at a cell centre. a. In the fourth quadrant, corners lie on the line $y=1-x$. \hskip 10 pt b. In the fourth quadrant, corners lie on the line $y=2-x$.
\end{figure}

The odd numbers one gets as areas of basic lattice units are only the sums of odd and even squares.
We see {\it for fabrics with rotational symmetry} why there are so few non-twill isonemal fabrics of small odd order, a fact noticed by Gr\"unbaum and Shephard in [3], p.~295, and given one explanation by them in [4], p.~158.
(The other contributing reason from [10] and [11] is that odd orders for non-twills with reflective symmetries must be products of relatively prime integers greater than 1, which obviously start at 15.)
Two sequences of examples of these lattice units with centre at a cell centre appear in Figure 4.
The sums of squares corresponding to the corners of the lattice units in Figure 4a appear down the main diagonal of the truncated infinite triangular array:
$$\left[\matrix{
4+1 & 9+4 & 16+9 & 25+16 & 36+25 & 49+36 & 64+49 & 81+64\cr
& 16+1 & 25+4 & 36+9 & 49+16 & 64+25 & 81+36 & 100+49\cr
&& 36+1 & 49+4 & 64+9 & 81+16 & 100+25 & 121+36\cr
&&& 64+1 & 81+4 & 100+9 & 121+16 & 144+25\cr
&&&& 100+1 & 121+4 & 144+9 & 169+16\cr
&&&&& 144+1 & 169+4 & 196+9\cr 
&&&&&& 196+1 & 225+4\cr 
&&&&&&& 256+1\cr} \right],\eqno{(2)}$$
whose $(i,j)$ entry (using matrix language) is missing if $i>j$ and is $M_1^2 +N_1^2$ if $i \leq j$, where $M_1 = i + j$ and $N_1 = j-i+1$.
The sums of squares corresponding to the corners of the lattice units in Figure 4b appear in the matrix down the superdiagonal.
Obviously the two smallest lattice units are the most oblique, the units with corners along each line $y=r-x$ tending asymptotically to conform to the strand directions as in the exceptional fabrics.
Likewise successive units with corners along horizontal lines tend toward having sides at $45^\circ$ to the strand directions, the slope of the lattice-unit sides in Roth types 1--32.
Note that every cell corner off the forbidden lines in the diagram is a corner of a basic lattice unit.
In \S 4 those that cannot give rise to isonemal prefabrics will be weeded out.

While placing one corner of the lattice units in the fourth quadrant and above the line $y=-x$ entailed no loss of geometric generality, reflecting those configurations in the forbidden lines is an obvious possibility.
The result of reflection in each of the forbidden lines is the same.
Since there is no reflective symmetry in the prefabrics being discussed here, this different placement of the lattice units would produce designs differing by such reflections from those we consider.
There is a handedness (cheirality) to these prefabrics that those with reflective symmetry (all Roth types 1--32) obviously lack.%
\footnote{There is a handedness to the {\it diagrams} in [10], but not to the designs.}
Beyond this section we shall consider lattice units of odd area, which are basic to all, so placed and so ignore (because we know they are there) all of the prefabrics that have symmetry groups based on their reflections.
We can think only about right gloves because we know so well what left gloves look like, or, if we prefer, we can think of these as left gloves and ignore the right; it makes no difference.
In the remainder of {\it this} section we cannot ignore this matter altogether.

Just as having corners of lattice units be at points of cell-centre/cell-corner type opposite to the centre ensures oddness of area, having them be at points of the same type $(k, \ell)$ with non-zero integers $k$, $\ell$, ensures evenness of area.
Staying off the forbidden lines is needed to ensure obliquity.
For definiteness and with no further loss of generality, consider $k > - \ell > 0$.
For these points, the lattice unit is the square on a hypotenuse running from
$(k, \ell)$ to its image under a quarter-turn about $(0,0)$, $(-\ell, k)$, the other sides of the triangle running from $(k, \ell)$ to 
$(k, k)$ and from there to $(-\ell,k)$ with lengths $k-\ell$ and $k + \ell$.
For any choice of parity of $k$ and $\ell$, the lengths $k+\ell$ and $k-\ell$ have equal parity, and so the area of the lattice unit, $(k-\ell)^2 + (k+\ell)^2$, is even.
Let us consider this quantity a little. 
It is equal to $2(k^2 + \ell^2)$, twice a sum of squares that can itself be odd or even.
If $k^2 + \ell^2$ is odd, then this unit is twice the area of a basic lattice unit.
If $k^2 + \ell^2$ is even, then this unit is twice the area of a lattice unit of one kind or the other of {\it these} lattice units, since an even sum of squares is the sum of two odd squares or of two even squares, and any pair of integers of the same parity is a case of $k \pm \ell$.
Halving within these units can obviously go on only so long; eventually a basic unit is reached.

I shall call a lattice unit $2^n$ times the area of a basic unit {\it of level} $n+1$, so that the basic level is level 1, as it is now provisionally defined (final definition in \S 4).
It will be handy to be able to refer to the dimensions of lattice units, and to do so in terms of the translations to which their sides correspond.
This notation is chosen to extend but conform to that of Roth [8].
Level 1 can be represented by the translation along its bottom edge $(M_1,N_1)$ of (1) with $M_1 > N_1$. 
The same lattice unit supplies many other translation vectors, of course, all linear combinations of $(M_1,N_1)$ and $(-N_1,M_1)$.
I have chosen the obliquity of the level-1 unit to be rotated a little anti-clockwise from the grid of the strands not just because $(M_1,N_1)$ is easier to write than $(M_1,-N_1)$ but also so that the fourth-quadrant corner positions in Figure 4 could match the corresponding number positions in (2) in an upper-triangular array.
Level 2 has analogous parameters $M_2 = M_1 + N_1$ and $N_2= M_1-N_1$, both odd.
Level 3 has
$$M_3 = M_2 + N_2 = 2M_1 \ \hbox{and}\ N_3= M_2-N_2 = 2N_1,\eqno(3)$$
both even.
Level 4 has $M_4 = 2M_1 + 2N_1 = 2M_2$ and $N_4= 2M_1-2N_1 = 2N_2$,
and level 5 has
$$M_5 = 2M_2 + 2N_2 = 2M_3 = 4M_1 \ \hbox{and}\ N_5= 2M_2-2N_2 = 2N_3 = 4N_1.\eqno(4)$$

Descriptions and drawings of these lattice units will be necessary and will follow some conventions.
Units of level 2 will be drawn escribing level-1 units the way the unit of Figure 2b escribes that of Figure 2a.
Units of level 3 will be drawn escribing units of level 2, and so on.
This relation is sometimes natural (as in Figure 2), and so I have made it conventional.
Note that lattice-unit corners at levels higher than 1 fall at cell centres rather than cell corners, which are used up by the corners of level-1 units.
In each quadrant, it puts the corners of units of odd level greater than 1 on the same side of the lines $y=\pm x$ as the corners of level-1 units in Figure 4a and corners of units of even level on the opposite side of the lines $y=\pm x$.
One nested set of units of five levels based on the smallest level-1 unit is illustrated in Figure 5a.

The basic translations, of which others are linear combinations, associated with the illustrated lattices at level 1 are $(M_1,N_1)$ and $(-N_1,M_1)$. 
Likewise $(M_3,N_3)$ and $(-N_3,M_3)$ and $(M_5,N_5)$ and $(-N_5,M_5)$.
But the different orientation of the even-level lattice units gives $(N_2,M_2)$, $(M_2,-N_2)$ and $(N_4,M_4)$, $(M_4,-N_4)$.
In Section 4 it will be necessary to see that level-1 lattice units as illustrated in Figure 5a can tessellate the plane in such a way as to place their corners on every strand boundary and that level-2 lattice units as illustrated place their corners in every strand.
The condition needed in both cases is the relative primality of $M_i$ and $N_i$, equivalently the existence of $m_i, n_i$, such that $m_iM_i + n_iN_i =1$.
To see the necessity, consider $(M_i, N_i)$ with a common factor $f$.
Lattice unit corners would then appear at most on every $f$th strand boundary or in every $f$th strand.
To see the sufficiency, the condition $m_iM_i + n_iN_i =1$ is a kind of recipe for how many translations $(M_1, N_1)$, $(-N_1, M_1)$, are required to move a lattice-unit corner one warp to the right, namely $m_1(M_1, N_1) - n_1(-N_1, M_1)$.
And likewise one weft up for level 1, $n_1(M_1, N_1) +m_1(-N_1, M_1)$.
For level 2, $m_2M_2 + n_2N_2 =1$ gives analogous recipes: 1 warp: $n_2(N_2, M_2) + m_2(M_2, -N_2)$; 1 weft: $m_2(N_2,M_2) -n_2(M_2,-N_2)$.
\begin{figure}
\begin{tabbing}

\noindent\epsffile{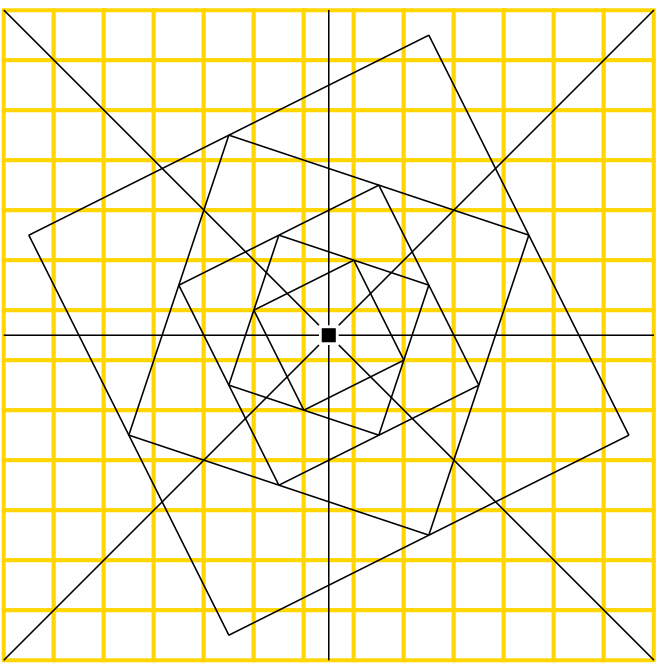}\hskip 5 pt \=\epsffile{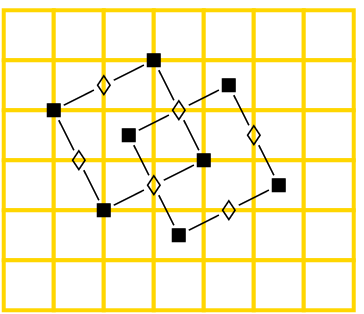}\\

\noindent (a)\>(b)
\end{tabbing}
\vskip -10 pt
\noindent Figure 5. a. A sequence of lattice units of levels 1 to 5 based on the smallest level-1 lattice unit of area 5.\hskip 10 pt b. Two single squares representing alternative tessellations of the plane by lattice units for the same transformation group.
\end{figure}

The discussion of the previous seven paragraphs is based entirely on placing the origin of the co-ordinate system at a cell centre.
If the origin is placed instead at a cell corner, all of the algebra carries over, and all points that are cell corners in Figures 4a, 4b, and 5a become cell centres and vice versa in Figure 5a.
Not only are there these two different ways to describe lattice units, but also they can both be used to describe some of the same weaving designs (e.g., that of Figure 7).
As Figure 5b illustrates, two level-1 lattices can describe the same transformation group, that is, place the same operations in the same places.
The placement of the lattice-unit sides, which are translations, are different, but the translations represented, $(M_1, N_1)$, $(-N_1, M_1)$, are the same.
I refer to this relationship between the two sets of tessellating squares as the {\it chain-mail} relation.
It almost appears that centring level-1 lattice units on both cell centres and cell corners is unnecessary.
That is not the case.
Some of the appearances are not deceiving, however.
\smallskip
\begin{Lem}{The mid-points of the sides of level-1 lattice units fall on the mid-points of cell sides.}
\end{Lem}
\noindent {\it Proof.} The side of each level-1 lattice unit is the hypotenuse of a right triangle with odd and even other sides.
Accordingly, its mid-point is an integral number of cell widths in the even direction and a non-integral number of cell widths in the odd direction from each end.
Whether beginning at a cell centre or cell corner, such a pair of displacements moves to the mid-point of a cell side, where `side' is being used in a general sense to include top and bottom.
\smallskip
\begin{Thm}{The centre of each level-1 lattice unit is of sort (cell centre or cell corner) opposite to its corners.}
\end{Thm}
\noindent {\it Proof.} Because the sides of the two sorts of level-1 lattice unit of the same size (as illustrated in Figure 5b) are each other's right bisectors, intersecting at a centre of half-turns according to Lemma 1, a path between corner and centre of either sort of lattice unit passes along half a side of each unit and therefore runs from a cell corner to a cell centre in one or the other direction.
\smallskip
\begin{Thm}{A square with corners at cell corners and side equal to the hypotenuse of a right triangle with other sides both odd and unequal or its reflection in a forbidden line is a lattice unit at level 2.}
\end{Thm}
\noindent {\it Proof.} Let the non-hypotenuse sides of the triangle be $2m+1$ and $2n+1$ in length for non-negative integers $m > n$.
Let the cell-corner ends of the hypotenuse be the origin and $(2m+1, 2n+1)$; then its mid-point is at $(m+{1\over 2}, n+{1\over 2})$, a cell centre.
The square of the hypotenuse is $4m^2+4m +4n^2 +4n +2 =2(2m^2 +2m +2n^2 +2n +1)$, where the second factor is obviously odd.
The line segment joining the mid-points of adjacent sides of the square on the hypotenuse is therefore the hypotenuse of some other right triangle whose other (vertical and horizontal) sides must be odd and even $m+n+1$ and $ m-n$. 
Their squares do sum to $2m^2 +2m +2n^2 +2n +1$.
The square on the smaller hypotenuse or its reflection in a forbidden line is therefore a level-1 lattice unit, and the original lattice unit is its escribing square.
The theorem is proved.

Let us consider lattice units with cell-corner corners since the two sorts of unit are related as they are.
Let a line segment from cell corner to cell corner be the side of a lattice unit and so be the hypotenuse of a right triangle with horizontal and vertical sides.
Its mid-point must be 

\noindent (i) at mid-side of a cell (other sides of the right triangle of which it is hypotenuse are even and odd) or

\noindent (ii) in the centre of a cell (other sides of the right triangle of which it is hypotenuse are both odd and unequal) or

\noindent (iii) at a corner of cells (other sides of the right triangle of which it is hypotenuse are both even and unequal) or

\noindent (iv) at a corner of cells or centre of a cell, but the other sides of the right triangle of which it is the hypotenuse are equal. This case is ruled out because the unit would not be oblique.

In case (i), the unit or its reflection in a forbidden line is at level 1.

In case (ii), the unit or its reflection in a forbidden line is at level 2 by Theorem 2.

In case (iii), the line segment is twice a cell-corner-to-cell-corner segment and so consider half of it, which can be of case (i) or (ii) or (iii) since it is oblique.
If it is of case (i), then the unit or its reflection in a forbidden line is of level 3, being four times the size of a level-1 unit.
If it is of case (ii), then the unit or its reflection in a forbidden line is of level 4, being four times the size of a level-2 unit.
If it is of case (iii), then it or its reflection in a forbidden line is of level 5 or higher and it can be halved again. 
The halving cannot go on indefinitely, and it stops in case (i) or (ii).
What we see is that any oblique $p4$ lattice unit or its reflection in a forbidden line is of level 1 or of a higher level based on a unit of level 1.
\begin{Thm}{Every oblique $p4$ lattice unit is of level 1 or of a higher level based on a level-1 lattice unit or the reflection of either in a forbidden line.}
\end{Thm}
\smallskip
As we shall see, levels higher than 5 are no use to us, but that is not because they cannot be the lattice units of periodic prefabric designs.
It is just that such prefabrics cannot be isonemal as I shall show in section 4 after showing in section 3 how the levels are used to distinguish the symmetry groups of isonemal prefabrics.
Because isonemality for order greater than 4 forces obliquity, we have proved the following theorem.
\begin{Thm}{A symmetry group of type $p4$ of an isonemal prefabric of order greater than 4 must have a lattice unit that is of level 1 or a higher level or is a reflection in a forbidden line of a lattice unit that is of level 1 or a higher level.}
\end{Thm}
\begin{figure}

\noindent\epsffile{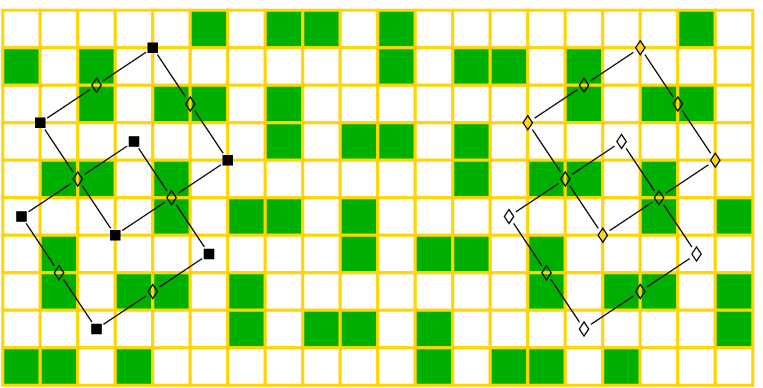}

\noindent Figure 6. The odd twillin 13-45-1 of species $36_1$ with $G_1$ marked on the left twice and $H_1$ marked correspondingly translated to the right by the order of the fabric.
\end{figure}
\vs
\noindent {\bf 3. Rotational symmetry groups and potential lattice units}

\noindent Because of the chain-mail way in which alternative level-1 lattice units for the same symmetry group can overlap, there are two alternative ways to describe the symmetry group illustrated in Figure 5b in terms of level-1 lattices.
One can place \blboxx s at the centres of cell-centred units and cell-corner \blboxx s at their corners, or one can place \blboxx s at the centres of point-centred units and cell-centre \blboxx s at their corners. 
Either way, one has placed \blbox at the same points spaced in exactly the same way because the corners and centre are of opposite cell-centre/cell-corner sorts.
The same configuration of centres of rotation can therefore be considered as having its level-1 lattice centred either on a cell or at cell corners as shown in Figure 5b.
One must also place \dia at the mid-points of the unit sides, which coincide. 
This placement in a level-1 lattice is one kind of group of the type that Roth numbers 36. 
It can be called the level-1 or odd-twillin subtype $36_1$.
It is exemplified by the square satin 5-1-1, to which there corresponds no twillin.
So to illustrate it with an odd twillin the next larger level-1 unit must be used, that of the next larger square satin, 13-1-1.
The illustration in Figure 6 is 13-45-1, marked on the left with $G_1$.
Its side-preserving subgroup $H_1$ is of type $p2$ and is marked on the right in Figure 6 in both of its versions.
These markings will serve as a description of the crystallographic type $p2$.
Its combination with $G_1$ of type $p4$, together with the presence of in-strand \diaa\hskip 4.5 pt (here some but not all \dia = \blboxx ${}^2$), is characteristic of isonemal prefabrics of species 36, of which two more subspecies will be defined at level 2.
Because fabrics of this type fall into genera I and III, the order is the same as the period, the common area of the $G_1$ and $H_1$ lattice units, 13 in the case of the example of Figure 6.

Technical information like that above is collected in Table 1 from where it is generated throughout this section.
The table is placed here so that its last column can be used as a navigational tool among the figures.
Unit areas in the table (columns 2 and 4) and orders (column 5) are multiples of the areas of the corresponding level-1 units.
It is the $H_1$ unit area that is the two-dimensional period.
\begin{table}
\begin{tabbing}
\hskip 40 pt \= $G_1$\hskip 15 pt \= $G_1$\hskip 18 pt\= $H_1$\hskip 15 pt\=\hskip 35 pt\=\\
Species                 \>unit \>unit \>unit \>Order \>Figures\\
                        \>area \>level \>area\\
$33_3$        		\>4    \>3    \>4        \>2  \>11b\\
$33_4$        		\>8    \>4    \>8        \>4  \>13b\\
34        		\>2    \>2    \>2        \>2  \>8a\\
$35_3$        		\>4    \>3    \>4        \>2  \>11a\\
$35_4$        		\>8    \>4    \>8        \>4  \>13a\\
$36_1$    		\>1    \>1    \>1        \>1  \>6\\
$36_2$    		\>2    \>2    \>2        \>2  \>8b\\
$36_s$			\>2    \>2    \>2        \>2  \>9\\
37        		\>8    \>4    \>16       \>4  \>14, 15\\
38        		\>4    \>3    \>8        \>4  \>12\\
39        		\>1    \>1    \>2        \>2  \>7, 16\\
\end{tabbing}
\vskip - 18 pt
\noindent Table 1. Relations among oblique lattice-unit sizes for species with rotational symmetry.
\end{table}

The type-$36_1$ placement of \blbox is not possible for \whbox because \whbox cannot be placed at the centre of a cell at all.
For \whbox within a cell to be a symmetry operation, that cell would have to be invariant under a quarter-turn (without $\tau$) but to keep its colour nevertheless, a plain contradiction.
(As indicated in [11], \diab is also not possible within a cell.)

\begin{figure}\noindent\epsffile{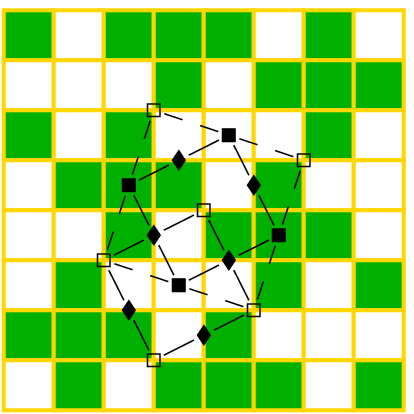}

\noindent Figure 7. The even twillin 10-93-1 of species 39 with two lattice units of $G_1$ outlined (level 1) and one lattice unit of $H_1$ outlined with dashes (level 2).
\end{figure}

The alternative ways of describing an allocation of operations are more evident when one says that placing \blbox at the centre of cell-centred units and cell-corner \whboxx s at their corners is the same as placing \whbox at the centre of corner-centred units and cell-centre \blboxx s at the unit's corners. 
For consistency, one must place \diabs (rather than \diaas \hskip -2.5 pt ) at the coincident mid-points of the unit sides.
This placement in a level-1 lattice is the kind of group action that Roth numbers 39.
It is illustrated in Figure 7, where the illustrated escribing square is the lattice unit for the side-preserving subgroup $H_1$, also of type $p4$.
Note that, while $G_1$ has the two lattices, the escribing squares of the one but not of the other are lattice units for $H_1$ because only one has \whbox at its centre as must be the case.
Roth illustrates species 39 with the only fabric of minimal order, 10-93-1.
Its order is the period, the area of the $H_1$ lattice unit, twice that of the $G_1$ unit.
This type is characterized by $G_1 \neq H_1$, both of type $p4$, with \diabs along strand boundaries in $G_1$ and without side-reversal $\hbox{\diaa} = \hbox{\blboxx\hskip 0.7 pt}^2$ in cells in $G_1$ and $H_1$.
Strands of this species have the most complicated symmetry possible; we return to it later.
This completes consideration of the two ways a level-1 lattice unit can be used.
The escribed squares of level-1 lattice units are the lattice units of the second level with side $d$ the diagonal of the level-1 lattice units. $M_1^2 + N_1^2 = d^2$.
Because they are escribed, they have twice the area of the level-1 unit, but each of their sides is the hypotenuse of a right triangle with other sides $M_1 + N_1$ and $M_1 - N_1$.
Both $M_1 + N_1$ and $M_1 - N_1$ are odd with $(M_1 + N_1)^2 + (M_1 - N_1)^2 = 2d^2$.
\begin{figure}
\begin{tabbing}

\noindent\epsffile{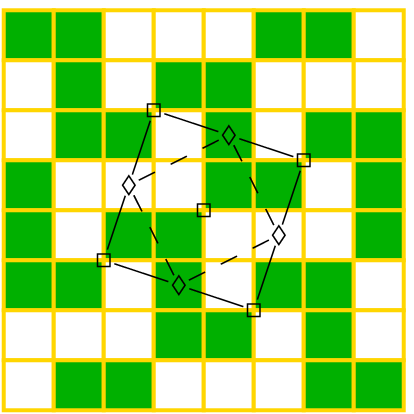}\hskip 10 pt \=\epsffile{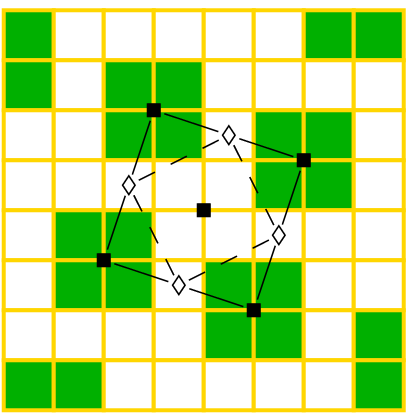}\\

\noindent (a)\>(b)\\
\end{tabbing}
\vskip - 18 pt
\noindent Figure 8. Two even twillins with lattice unit for $G_1$ outlined (level 2) and the level-1 unit on which it is based dashed. a. 10-107-1 of species 34. \hskip 10 pt b. 10-27-1 of species $36_2$.
\end{figure}

We turn to second-level units, first with centre at cell corners, accordingly with corners also at cell corners.
The mid-points of the units' sides fall at the centres of cells and so can support only \dia (no $\tau$ as when \whbox and \blbox are both used as in species 39).
This limitation requires \whbox at both centre and corners or \blbox in both places.
Using \whbox produces a symmetry of Roth type 34 with $H_1 = G_1$ illustrated by 10-107-1 of the smallest possible order (Figure 8a).
This type is characterized by the equality of $H_1$ and $G_1$ and the presence of half-turn centres at the middle of the units' sides in cells.
Using \blbox instead produces another symmetry of Roth type 36, which he illustrates with 10-27-1 (Figure 8b).
This species, which I number $36_2$, is clearly distinct from species $36_1$, whose lattice unit is of level 1.
The side-preserving subgroup, which is of type $p2$, is the same restriction of the operations \blbox to $\hbox{\blboxx\hskip 0.7 pt}^2 = \hbox{\dia}$ as in Figure 6b, the lattice unit for $H_1$ being that for $G_1$.
As both of these species are twillins, their order equals the period, the area of the common lattice unit.
\begin{figure}

\noindent\epsffile{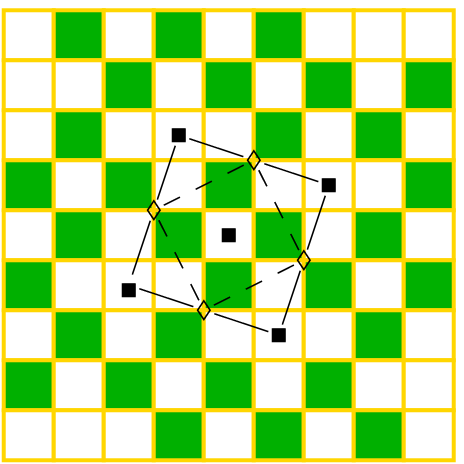}

\noindent Figure 9. The even twillin 10-85-1 of species $36_s$ with a lattice unit of $G_1$ outlined (level 2) and the level-1 unit on which it is based dashed.
\end{figure}
Second-level units with centre and corners at cell centres can support only \blbox there, and this gives groups of a third subtype $36_s$ of Roth type 36.
Since the other subscripts are chosen to reflect level, s is chosen for several reasons, one of which is that it looks more like a 2 than like a 1, 3, or 4.
The mid-points of sides must be equipped with \diaas\hskip -2.5pt , which fall at cell corners.
This species is illustrated by 10-85-1 in Figure 9, in which the central but not the corner pale cells are surrounded by dark crosses.
The existence of the three subtypes 1, 2 and s was remarked on by Roth [7], where he distinguishes among them by the fact that the strands' half-turn symmetries (\diaa\hskip 1.6 pt) arise as $\hbox{\blboxx\hskip 0.7 pt}^2$ in some, none, and all cases, respectively.
One can see in Figures 6, 8b, and 9 that the distinction is also a matter of where the \blboxs and mid-side \diaas fall.
Twillins of odd order are of species $36_1$, and twillins of even order fall into `none' or `all' depending on whether their lattice unit is corner-centred or cell-centred.
\smallskip

\begin{Thm}{Square satins of odd and even order fall into species $36_1$ and $36_s$ respectively.}
\end{Thm}

This concludes what can be done with lattice units at the second level.
The above also concludes the use of cell-centred lattice units, which {\it must} be used only for species $36_s$, another reason for its distinctive subscript. And s is for satin.
\begin{figure}

\noindent\epsffile{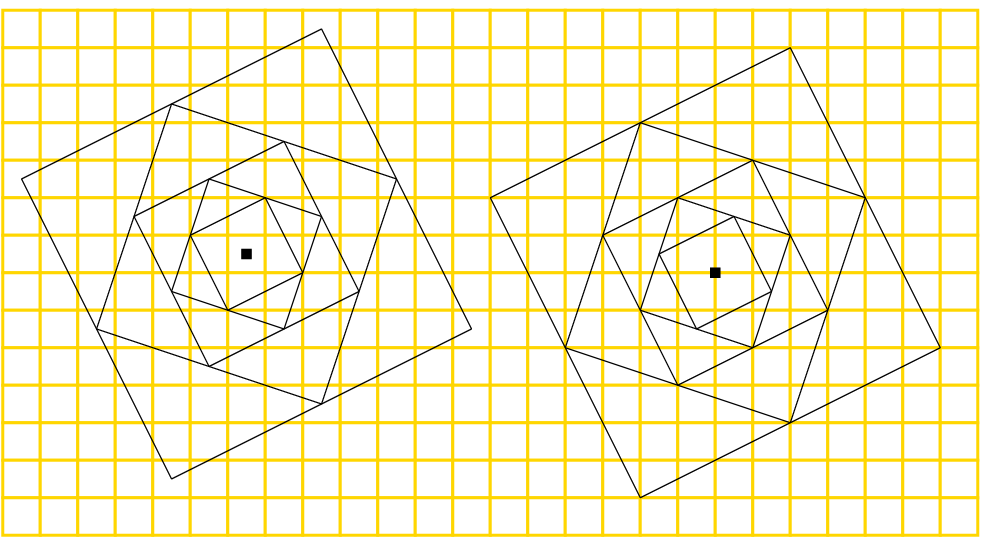}

\noindent Figure 10. Potential lattice units of levels 1 to 5 based on the smallest level-1 unit, cell-centred on the left and point-centred on the right.
\end{figure}
\begin{Thm}{Lattice units of symmetry groups of isonemal prefabrics at level 3 and higher have centres and corners at cell corners.}
\end{Thm}
\noindent {\it Proof.} Lattice units of levels higher than one have corners {\it and} centres either both at cell centres or both at cell corners.
At level 2 this does not prevent their application to isonemal prefabrics because the centres of half-turns at mid-side suffice to relate strands to adjacent strands.
But at level 3, the mid-side half-turn centres (at the corners of the inscribed level-2 units) are also at cell centres, leaving nothing to relate adjacent strands.
Translations cannot do it, because the distances involved, $M_3$ and $N_3$, are even by (3).
Warps and wefts are interchanged, but the prefabrics are not isonemal.
This failing is inherited by all higher levels of potential lattice units with cell-centre centres.
The theorem is proved.
\begin{table}
\begin{tabbing}
\hskip 0.1 in \= Centre \hskip 0.4 in \= Corners \hskip 0.55 in \= Mid-sides \hskip 0.38 in \= Case\\
\smallskip
\>\blbox \>\blbox \>\dia \> (\blboxx , \blboxx )\\

\>\whbox \>\whbox \>\dia \> (\whboxx , \whboxx )\\

\>\whbox \>\blbox \>\diab \> (\whboxx , \blboxx )\\
\end{tabbing}
\vskip - 18 pt
\noindent Table 2. Rotation assignments to lattice-unit positions for levels 3 and 4.
\end{table}

We move on to consider lattice units of level 3 with cell-corner centres.
Because the corners and mid-sides of these units also fall at cell corners, dimensions all being even, there is complete freedom of assignment of \whbox and \blbox to the corners subject only to consistency. 
The configurations mixing \whbox and \blbox are merely alternate ways of specifying the same array of operations in terms of different lattice units arranged like chain mail.
The possibilities are displayed in Table 2.
We treat the alternative with \whbox at the centre and \blbox at the corners so that the lattice unit for $H_1$ can escribe that for $G_1$.
\begin{figure}
\begin{tabbing}
\noindent\epsffile{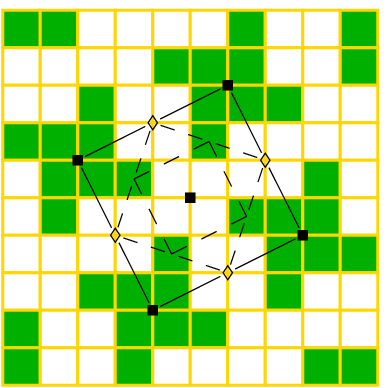}\hskip 5 pt\=\epsffile{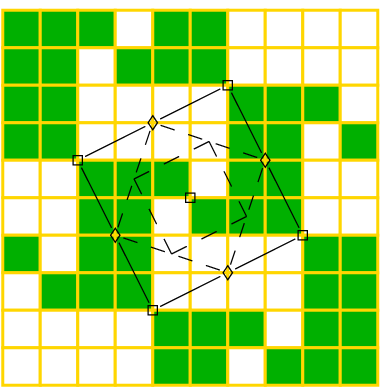}\\

\noindent (a)\>(b)\\
\end{tabbing}
\vskip -18 pt
\noindent Figure 11. Genus-III examples with $G_1$ outlined (level 3) and the level-1 unit and intermediate level-2 unit dashed. a. 10-39-1 of species $35_3$.\hskip 10 pt b. 10-55-2 of species $33_3$.
\end{figure}

The (\blboxx , \blboxx ) case has side-preserving subgroup expressed analogously as ($\hbox{\blboxx\hskip 0.7 pt}^2$, $\hbox{\blboxx\hskip 0.7 pt}^2$) with the same lattice unit, so that the $G_1/H_1$ type pair is $p4/p2$.
As all rotation operations fall at cell corners, the strands have only the trivial (translation by order) symmetry group, and so the Roth type is 35, illustrated by Roth with 10-39-1, which is as small an order as possible (Figure 11a).
This example is the first fabric here that is neither a satin nor a twillin.
All designs with this or larger lattice units will be of genus III alone, which is the case here, or III, IV, and V together, or V alone because the centres of rotation are spaced out sufficiently far to create what are in effect two classes of strands alternating.
Because $H_1$ has the same lattice unit as $G_1$, which is of level 3, the period is four times that of the level-1 unit, and because the genus is III the period is twice the order, which is therefore twice the area of the level-1 unit.
The operations are exactly twice as far apart as those of the level-1 lattice unit, which has been dashed---along with the intermediate level-2 unit---in Figure 11a.
Because it is based on a level-three lattice unit, this species will be called $35_3$.

$G_1$ in the (\whboxx , \whboxx ) case is its own side-preserving subgroup.
Again the strand symmetry group is trivial (only translation by order), and so the Roth type is 33, which Roth illustrates with 10-55-2, which is as small as possible.
In Figure 11b illustrating it, the level-2 and level-1 lattice units have been dashed.
Because it is based on a level-three lattice the species in this case will be called $33_3$.
All spacing being the same as in the (\blboxx , \blboxx ) case, the order and period are the same as there.

The (\whboxx , \blboxx ) case has the same operations as Roth type 39 (Figure 7) but on a larger scale.
As in type 39, its side-preserving subgroup cannot use the \diabs, and so the half-turns at mid-side of the $H_1$ lattice unit must be squares of the corner \blboxx s, making the $H_1$ lattice unit's corners be the \whbox centres of four (level-3) $G_1$ lattice units surrounding the $G_1$ unit.
Accordingly, the subgroup's lattice unit is at level 4.
As we shall see, the minimum order of such a design is 20, and Roth had to create a design of this size to illustrate this species, 38.
Figure 12a reproduces the design of his [7], Figure 8,
with the $G_1$ lattice unit as well as the smaller lattice units shown.
The $H_1$ lattice unit (with $\hbox{\blboxx\hskip 0.7 pt}^2 = \hbox{\dia}$) is shown in Figure 12b.
As in type 39, $G_1 \neq H_1$, both of type $p4$, making the symbol $p4/p4$.
Because the genus is III, the period is twice the order, but because the $H_1$ lattice unit is at level 4 the period is 8 times the level-1 unit, making the order 4 times the level-1 unit, the minimum being 20.
The type can be distinguished from 39 by its lattice-unit level, but this is not how the types were distinguished by Roth, who used a consequence.
As Figure 7 indicates, because the lattice unit of type 39 has \diab within strands, the strand symmetry that Roth calls 12/12 is imposed on strands.
With the larger lattice unit of type 38, \diab being at cell corners, the strand symmetry is weakened to what he calls 11/11, translation with side reversal ($\tau$) instead of half-turn with side reversal ($\tau$).
This translation is the result of $\hbox{\dia} = \hbox{\whboxx\hskip 0.7 pt}^2 = \hbox{\blboxx\hskip 0.7 pt}^2$ and \diab on the same edge of each strand alternating at intervals of order/4 visible in Figure 12a (twice vertically and twice horizontally, all with $\hbox{\dia} = \hbox{\blboxx\hskip 0.7 pt}^2$).
We have now exhausted the potential of the third level.
\begin{figure}
\begin{tabbing}
\noindent\epsffile{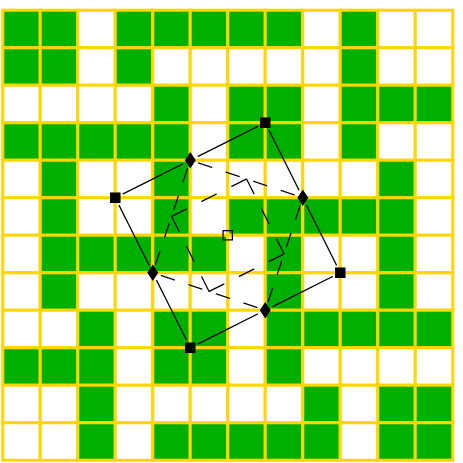}\hskip 5 pt \=\epsffile{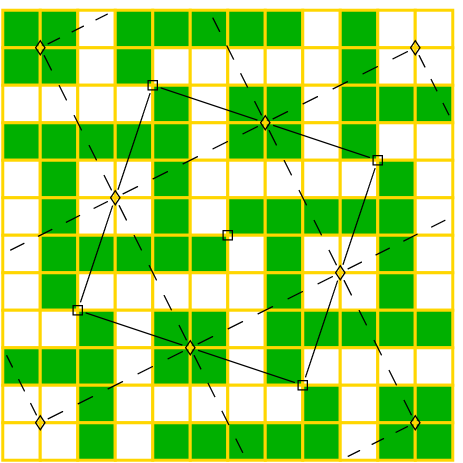}\\

\noindent (a)\>(b)\\
\end{tabbing}
\vskip -18 pt
\noindent Figure 12. Roth's example 20-19437 of species 38, after [7], Figure 8. a. $G_1$ (level-3) lattice unit outlined with the level-1 and intermediate (level-2) units dashed.\hskip 10 pt b. $H_1$ (level-4) lattice unit outlined with $G_1$ lattice units dashed.
\end{figure}

The level-4 groups have the same freedom of operation choice as those at level 3 (Table 2).
As at level 3, the (\blboxx , \blboxx ) case has $G_1/H_1$ type pair $p4/p2$ and is in other ways entirely similar to subtype $35_3$; it just has a lattice unit twice as large.
Its minimum possible order is accordingly 20.
Figure 13a illustrates the corresponding species, which I call $35_4$, with a design that is 10-85-1 doubled (cf.~Figure 9).
Though it is less easy to see, just as the design of Figure 9 without the centres of the dark crosses pale would be that of Figure 7 (10-93-1), so that of Figure 13a without the centres of the dark crosses pale would have the $G_1$ illustrated in Figure 12a in its alternative representation with \blbox at unit centres (and the \diaas of Figure 13a changed to \whboxx s at unit corners).
\begin{figure}
\begin{tabbing}
\noindent\epsffile{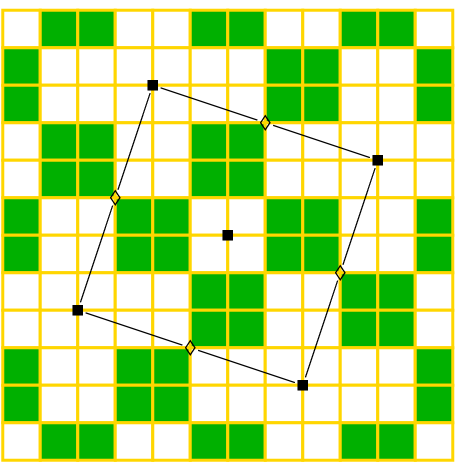}\hskip 5 pt\=\epsffile{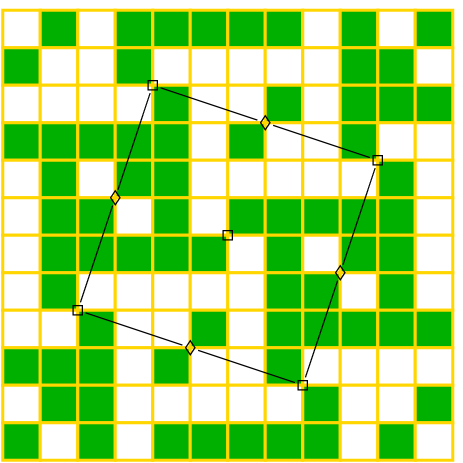}\\

\noindent (a)\>(b)\\
\end{tabbing}
\vskip -18 pt
\noindent Figure 13. $G_1$ lattice units for level-4 fabrics of order 20. a. 10-85-1 doubled of species $35_4$. b. Example of species $33_4$.
\end{figure}

The (\whboxx , \whboxx ) case has a group of the same type $p4/-$ as that in Figure 11b but at the higher level.
It is its own side-preserving subgroup, and the centres of rotation are sufficiently spaced out to have no implications for strand symmetry.
Accordingly its Roth type is 33.
Because of its level-4 lattice unit, I distinguish the species as being $33_4$.
It is illustrated in Figure 13b by a design just different enough from that of Figure 12a to have its symmetry group be this particular subgroup of $G_1$ of the design of Figure 12a.
This group is also a subgroup of that of the design of Figure 11b of subtype $33_3$. 
The minimum order for this species is 20 because its lattice unit is twice the size of that for subtype $33_3$.
The subdivision into species $33_3$ and $33_4$, $35_3$ and $35_4$, is the refinement of Roth's taxonomy in this paper.

When we turn to the (\whboxx , \blboxx ) case at level 4, the operations match Roth types 39 and 38 (Figures 7 and 12 respectively), and so $G_1 \neq H_1$, both of type $p4$.
But now the centres of rotation are sufficiently spaced out that the symmetry group of the strands is trivial (translation only).
These are characteristics of Roth's type 37.
The design that Roth invented to illustrate this species, 20-3391, is illustrated in Figure 14.
Its side-preserving subgroup has, like types 39 and 38, a lattice unit of the next level, in this case level 5.
It is illustrated for another design of this type in Figure 15.
The exceptional prefabric 4-1-2* is the smallest (and only small) example of this Roth type (Figure 2).
They fall into the same type because the centres of rotation are spaced out sufficiently that they have no implications for symmetry of the strands beyond periodicity.
The obliqueness of these larger lattice units allows them to be as large as we please, unlike non-oblique units.
We have exhausted the possibilities of lattice units of level 4, which is all the possibilities there are, as the next section will prove.
\begin{figure}
\noindent\epsffile{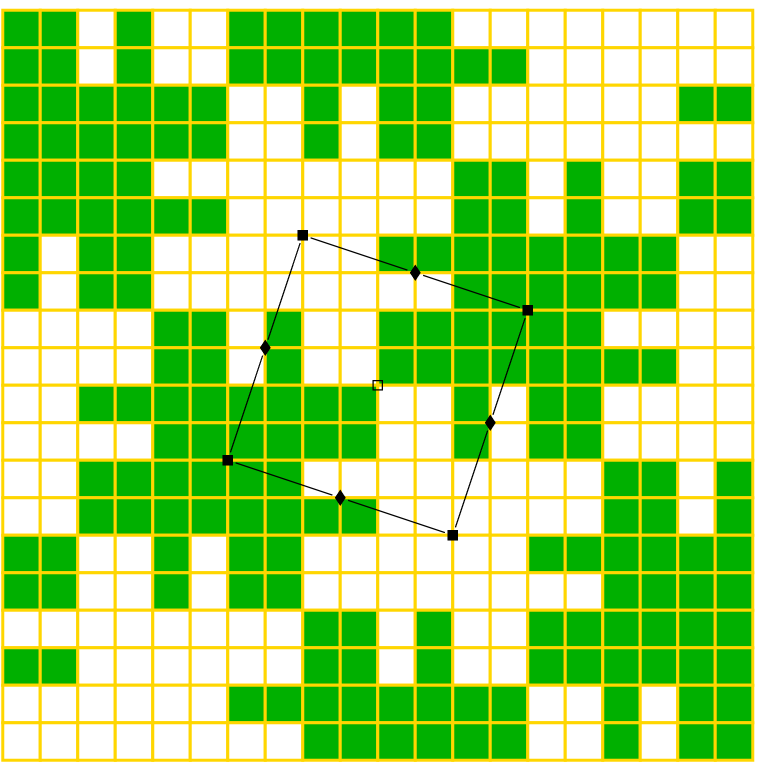}

\noindent Figure 14. The genus-V order-20 example of species 37 invented by Roth with $G_1$ (level 4) outlined. After [7], Figure 7.
\end{figure}
\begin{figure}
\noindent\epsffile{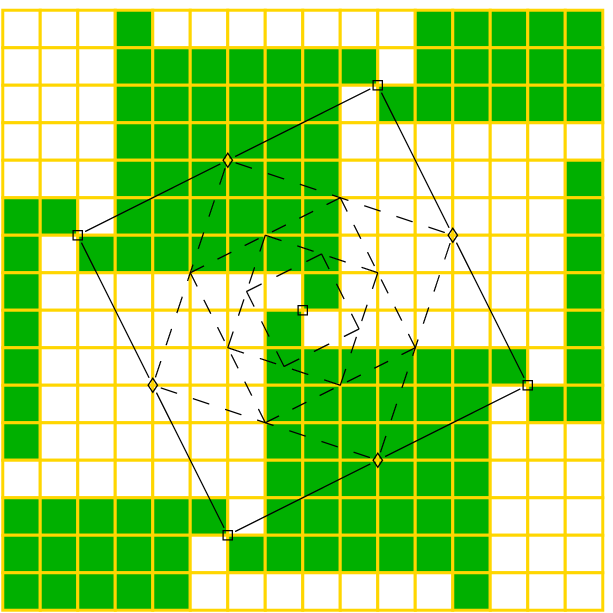}

\noindent Figure 15. A second order-20 example of species 37 with $H_1$ of level 5 marked and units of level 1 to 4 dashed.
\end{figure}

\noindent {\bf 4. Limitations on Lattice Units}

\noindent Limitations on usefulness of lattice units are of two kinds, which level-1 units work, and which levels built on them work.
As we have already seen that levels above the second built on level-1 units centred on a cell cannot be used, we initially pursue this sort of consideration.
As the lattice-unit level rises, the centres of rotation are spaced out more and more, with the result that at level 5 the group is no longer able to serve as the symmetry group of an isonemal prefabric.
As in Figure 2b, so in Figure 15 (where it is the side-preserving subgroup $H_1$ whose level-5 lattice unit is not dashed), pairs of adjacent strands not containing \whbox between them are transformed into pairs of strands by $H_1$, but there is no operation interchanging the members of those pairs (the job of the \diabspecial of $G_1$ as illustrated in Figure 14).
That this is always the case is because the distance along the edge of a unit at level 5 from corner to mid-side is twice the whole edge length of the unit at level 1.
Overlapping level-1 units are illustrated in Figure 5b. 
In that figure, one sees that level-1 lattice-unit corners at cell corners suffice to transform wefts to adjacent wefts (as do the half-turns at mid-sides).
When the side lengths are doubled (level 3), the strands fall into two classes; this was the level at which genera I and II ceased to appear as noted at the time.
With another doubling, the two classes cease to be transformed into each other.
A prefabric with such a group would not be isonemal.
While this is plain enough for level-1 lattice units as illustrated in Figure 5b, where the right triangle of which the lattice unit is the square on the hypotenuse has other sides 1 and 2, it is not obvious when the square on the hypotenuse is not on the main diagonal of array (2).
A proof, accordingly, will have to wait until we have studied the entries of (2) further (Lemma 9).

To see which of the entries of array (2) are useful and how their usefulness at level 1 fails to extend to level 5---although it does extend to level 4---we must see how numbers between row one and the main diagonal behave.
The displayed array (2) is the beginning of an infinite array of all odd sums of squares $M_1^2 +N_1^2$.
Those that cannot be used are those where $M_1^2$ and $N_1^2$ are not relatively prime, equivalently $M_1$ and $N_1$ are not relatively prime as noted two paragraphs before Lemma 1.
The useful necessary and sufficient condition for relative primality of $M_1$ and $N_1$ is the existence of integers $p$ and $q$ such that $pM_1 + qN_1 =1$.
Along the first row of (2), where $M_1=N_1+1$, this condition is obviously satisfied with $p=1$ and $q=-1$. 
Equally easily, in the second row $36+9$, $81+36$, $144+81$, and every third sum has 3 as common divisor of $M_1$ and $N_1$.
This is a special case of a general phenomenon, as I now show.
As specified at (2), in the general position $(i,j)$, $M_1=i+j$ and $N_1=j-i+1$.
Along the $i$th row, $M_1$ and $N_1$ begin with a difference of $2i-1$ and each of $M_1$ and $N_1$ increases by 1, column by column.
Accordingly, for $i>1$, every $(2i-1)$th entry, $M_1$ and $N_1$ have $2i-1$ as a common factor, starting with the $(2i-1)$th.
The $(2i-1)$th entry in the $i$th row is in position $(i,3i-2)$ and so has
$M_1=4i-2$ and $N_1= 2i-1$, and as $j$ is increased by $2i-1$, both $M_1$ and $N_1$ increase by that amount, remaining multiples of $2i-1$. 
But they also have the prime factors $\ell$ of $2i-1$ as factors every $\ell$th entry so that row 8, for example, with $2i - 1= 15$ loses every third {\it and} every fifth entry.
The entries behave as though the absent subdiagonal had entries $M_1=2i-1$, $N_1=0$.
If $2i-1 = m\ell$, then, for every $k$, $k\ell$ positions across among the actual entries
$M_1=(m+k)\ell$ and $N_1=k\ell$ have $\ell$ as a common factor.
The array with only its useful numbers remaining begins as follows:
$$\left[\matrix{
4+1 & 9+4 & 16+9 & 25+16 & 36+25 & 49+36 & 64+49 & 81+64\cr
& 16+1 & 25+4 &  & 49+16 & 64+25 &  & 100+49\cr
&& 36+1 & 49+4 & 64+9 & 81+16 & & 121+36\cr
&&& 64+1 & 81+4 & 100+9 & 121+16 & 144+25\cr
&&&& 100+1 & 121+4 &  & 169+16\cr
&&&&& 144+1 & 169+4 & 196+9\cr 
&&&&&& 196+1 & 225+4\cr 
&&&&&&& 256+1\cr}\right].\eqno{(5)}$$
It is pairs of terms (representing lattice-unit orientations) and not just sums (mere lattice-unit areas) that are rejected: while 125 as $100+25$ must be rejected, 125 as $121+4$ is acceptable.
There follows some discussion of how relative primality relates to isonemality with the goal of proving Theorem 7 in a couple of pages.

\begin{Lem}{If an isonemal fabric of species $36_1$ or $39$ has a lattice unit centred on a cell centre (and so with cell-corner corners), it is of level 1 with $M_1$ and $N_1$ relatively prime, and conversely if a fabric has symmetry group with lattice unit that is a square centred on a cell centre and that is on the hypotenuse of a right triangle with horizontal and vertical sides that are relatively prime and of opposite parity, $M_1$ and $N_1$, then it is isonemal and of species $36_1$ or $39$.}
\end{Lem}
{\it Proof}.
The first half was proved before Lemma 1, and so we turn to the second.
The opposite parities of $M_1$ and $N_1$ make the lattice unit be of level 1, and their relative primality, through $pM_1 + qN_1=1$, makes the central cell of the lattice unit translatable to every strand, vertical or horizontal. 
The fabric is accordingly isonemal.
It is of species $36_1$ if each corner of the lattice unit has \blboxx , and of species 39 if each corner has \whboxx .

\begin{Lem}{If an isonemal fabric of species $36_1$ or $39$ has a lattice unit centred on cell corners (and so with cell-centre corners), it is of level 1 with $M_1$ and $N_1$ relatively prime, and conversely if a fabric has symmetry group with lattice unit that is a square centred on cell corners and that is on the hypotenuse of a right triangle with horizontal and vertical sides that are relatively prime and of opposite parity, $M_1$ and $N_1$, then it is isonemal and of species $36_1$ or $39$.}
\end{Lem}
{\it Proof}.
The first half was proved before Lemma 1, and so we turn to the second.
The opposite parities of $M_1$ and $N_1$ make the lattice unit be of level 1, and their relative primality, through $pM_1 + qN_1=1$, makes each corner cell of the lattice unit translatable to every strand, vertical or horizontal. 
The fabric is accordingly isonemal.
It is of species $36_1$ if the centre of the lattice unit has \blboxx , and of species 39 if the centre has \whboxx .

\smallskip
Lemmas 2 and 3 describe different ways of speaking of the same fabrics.
\begin{Lem}{The (odd) sum and difference of relatively prime numbers of opposite parity are always relatively prime, and conversely, if the sum and difference of numbers are relatively prime, then those numbers themselves are relatively prime without any extra parity assumption.}
\end{Lem}
{\it Proof.} Let relatively prime even and odd numbers be $p$ and $i$ (for the French).
Then there exist coefficients $c_p$ and $c_i$ such that $c_p p +c_i i = 1$.
While $c_i$ must be odd, $c_p$ can be either odd or even.
If $c_p$ is even, then $c_p^\prime = c_p + i$ and $c_i^\prime = c_i -p$ are different coefficients that are both odd and make $c_p^\prime p +c_i^\prime i = 1$.
So without loss of generality, let $c_p$ and $c_i$ both be odd.
Then $(c_p \pm c_i)/2$ are both integers as would not be the case with $c_p$ and $c_i$ of opposite parities.
But the identity
$${c_p + c_i \over 2}(p+i) + {c_p - c_i \over 2}(p-i) = c_p p +c_i i = 1$$
shows that $p\pm i$ are relatively prime.

It is obvious that if $c(a+b) + d(a-b) = 1$, then $(c+d)a + (c-d)b = 1$, and difference of parity of $a$ and $b$ is built into the assumption.
\begin{Cor}{In the present context there is equivalence between the relative primality of $M_1$ and $N_1$ and of $M_2$ and $N_2$.}
\end{Cor}

\begin{Lem}{An isonemal fabric of species $34$ or $36_2$ has a level-2 lattice unit centred on cell corners (and so with cell-corner corners) with $M_2$ and $N_2$ relatively prime, and conversely if a fabric has symmetry group with lattice unit that is a square centred on cell corners and that is on the hypotenuse of a right triangle with horizontal and vertical sides that are relatively prime and odd, $M_2$ and $N_2$, then it is isonemal and of species $34$ or $36_2$.}
\end{Lem}
{\it Proof.}
A level-2 unit centred on cell corners has centres of half-turns at mid-side cell centres configured as the corners of a level-1 unit.
Those in-cell centres are distributed to every strand if and only if there are $p$ and $q$ such that $pM_1 + qN_1 = 1$.
The symmetry group is transitive on strands only if $M_1$ and $N_1$, and equivalently $M_2$ and $N_2$ by Lemma 4, are relatively prime. 
Conversely, if the sides of the right triangle are odd and relatively prime, then the square with corners at cell-centre mid-sides is centred at cell corners and is of level 1 with $M_1$ and $N_1$ relatively prime.
The relative primality of $M_1$ and $N_1$ assures the distribution of its corner half-turn centres to every strand, which assures isonemality.
If the quarter turns in such a configuration are \whboxx s then the species is 34, and if they are \blboxx s then the species is $36_2$.
\begin{Lem}{An isonemal fabric of species $36_s$ has a level-2 lattice unit centred on a cell centre (and so with cell-centre corners) with $M_2$ and $N_2$ relatively prime, and conversely if a fabric has symmetry group with lattice unit that is a square centred on a cell centre and that is on the hypotenuse of a right triangle with horizontal and vertical sides that are relatively prime and odd, $M_2$ and $N_2$, then it is isonemal and of species $36_s$.}
\end{Lem}
{\it Proof.} 
A level-2 unit centred on a cell centre has centres of half-turns at mid-side cell corners configured as the corners of a level-1 unit.
Those half-turn centres, which do all of the strand-to-adjacent-strand rotations, are distributed to every strand boundary if and only if there are $p$ and $q$ such that $pM_1 + qN_1 = 1$.
Every strand is accordingly related to adjacent strands by half-turns only if $M_1$ and $N_1$, and equivalently $M_2$ and $N_2$ by Lemma 4, are relatively prime. 
Conversely, if the sides of the right triangle are odd and relatively prime, then the square with corners at cell-corner mid-sides is centred at a cell centre and is of level 1 with $M_1$ and $N_1$ relatively prime.
The relative primality of $M_1$ and $N_1$ assures the distribution of its corner half-turn centres to every strand boundary, which assures isonemality.
The species of a fabric with symmetry group so configured is $36_s$.
\begin{Lem}{An isonemal fabric of species $33_3$, $35_3$, or $38$ has a level-3 lattice unit based on a level-1 lattice unit with $M_1$ and $N_1$ relatively prime, and conversely, a fabric with symmetry group having lattice unit of level 3 based on a level-1 lattice unit with $M_1$ and $N_1$ relatively prime is isonemal and of species $33_3$, $35_3$, or $38$.}
\end{Lem}
{\it Proof.}
The level-1 lattice unit within the level-2 unit within the level-3 unit is cell-corner centred. 
A congruent level-1 square with cell centre is one quarter of the level-3 lattice unit.
Since the corners, mid-sides, and centres of the level-3 lattice unit (half-turn centres or quarter-turn centres with or without $\tau$) regarded as centres of half-turns for this purpose, lie at the corners of the level-1 squares with cell-centres and tessellating the plane with corners on every strand boundary, $M_1$ and $N_1$ must be relatively prime.
Conversely, if $M_1$ and $N_1$ are relatively prime, then those level-1 squares tessellate the plane with corners on every strand boundary, and isonemality follows.
The fabric is then of species $33_3$, $35_3$, or 38 accordingly as the quarter-turns are all \whboxx s, all \blboxx s, or a mixture of \whboxx s and \blboxx s.
\begin{Lem}{An isonemal prefabric of species $33_4$, $35_4$, or $37$ has a level-4 lattice unit based on a level-1 lattice unit with cell-corner centre and $M_1$ and $N_1$ relatively prime, and conversely, a prefabric with symmetry group having lattice unit of level 4 based on a level-1 lattice unit with cell-corner centre and $M_1$ and $N_1$ relatively prime is isonemal and of species $33_4$, $35_4$, or $37$.}
\end{Lem}
{\it Proof.}
Consider one level-4 lattice unit of a design of species $33_4$, $35_4$, or $37$.
It has cell-corner centre, corners, and mid-sides, all of them being the corners of the level-2 lattice unit that is one quarter of it.
For the tessellation of the plane to distribute corners, half-turn centres and quarter-turn centres with or without $\tau$ (regarded as centres of half-turns for this purpose) to all strand boundaries, $M_2$ and $N_2$ must be relatively prime.
But the relative primality of $M_2$ and $N_2$ is equivalent to that of $M_1$ and $N_1$ by Lemma 4, and so the corresponding level-1 lattice unit has $M_1$ and $N_1$ relatively prime.
By the same equivalence, the relative primality of $M_1$ and $N_1$ ensures that what can be regarded as half-turn centres are suitably distributed, as the corners of the level-2 lattice units, to all strand boundaries so that the prefabric is isonemal.
It is then of species $33_4$, $35_4$, or 37 accordingly as the quarter-turns are all \whboxx s, all \blboxx s, or a mixture of \whboxx s and \blboxx s.

\begin{Lem}{A symmetry group with lattice unit at level 5 or higher with cell-corner centre cannot be the group of an isonemal prefabric.}
\end{Lem}
{\it Proof.} As examples have already shown, a unit at level 3 or 4 with cell-corner centre {\it can} be the lattice unit of an isonemal prefabric.
The geometry is that of Theorem 6 but for pairs of strands instead of individual strands.
As in Figure 10b, consider adjacent pairs of strands with the centre, mid-sides, and corners of the level-5 lattice unit between the strands of the pairs.
This is possible because half-sides have components $M_3$ and $N_3$ even.
One can see then that the configuration of Figure 10a at level 3 for single strands is duplicated in Figure 10b at level 5 for those adjacent pairs of strands.
All centres of rotation lie within the pairs, and nothing relates adjacent pairs to each other.
The translations along the sides of the lattice unit are too large, being multiples of four in their horizontal and vertical components, to move any pair to an adjacent pair.
This is and must be the case because the distances involved at level 5, $M_5$ and $N_5$, are twice those at level 3 (themselves even) by (4).
The lemma is proved.

\smallskip
Since lattice units of level 1 are of use only if they have $M_1$ and $N_1$ relatively prime as well as being of opposite parity, it is reasonable to build relative primality of $M_1$ and $N_1$ into the definition of {\it level one} and also of higher levels.
\begin{Def} An oblique $p4$ lattice unit will be called of {\em level one} if it is the square on the hypotenuse of a triangle with horizontal and vertical sides odd and even and relatively prime.
\end{Def}

Collectively Theorem 6 and Lemmas 2, 3, and 5--9 give the following theorem.
\begin{Thm}{The symmetry group of type $p4/-, p4/p4$, or $p4/p2$ of an isonemal prefabric must have lattice unit either with cell-centre centre and of level 1 or 2 or with cell-corner centre and of level 1, 2, 3, or 4.}
\end{Thm}
While this analysis involves only a slight refinement of the taxonomy of Roth, some complexity has been added in considering groups rather than types of groups, and that needs to pay its way.
It does this primarily in the geometrical insight that it offers.
But it is also useful in consideration of questions about weaving patterns.
In the following three sections, the questions taken up in [10] and [11] will be considered for these species.
Another question altogether is discussed in Section 8.

In view of the completion of the project of extending the work of Richard Roth on the symmetry groups of isonemal fabrics, it may be useful to present here an extension of the navigational part of Table 1 backwards through species 1--32. This is done in Table 3.
\begin{table}
\begin{tabbing}
Species\hskip 10 pt \= Figures in [10]\hskip 10 pt \= Species\hskip 10 pt\= Figures in [11]\\
$1_e$        	\>10a             \>11    \>5a\\
$1_m$        	\>4a              \>12    \>5b, 6\\
$1_o$        	\>4b, 10b         \>13    \>7a\\
$2_e$        	\>12b             \>14    \>7b\\
$2_m$    	\>12a             \>15    \>8\\
3    	        \>8a              \>17    \>9a\\
$4_e$    	\>8b              \>18    \>9b\\
$4_o$		\>13              \>19    \>9c\\
$5_e$        	\>2a, 14a         \>21    \>10\\
$5_o$        	\>14b, c          \>22    \>11a\\
6    	 	\>6               \>23    \>11b\\
$7_e$        	\>16a             \>25    \>2a\\
$7_o$        	\>15              \>26    \>2b, 4\\
$8_e$        	\>9a, 17a         \>27    \>12\\
$8_o$        	\>9b, 16b, 17b    \>28    \>13, 14\\
9    	        \>7               \>29    \>15a\\
10    		\>5               \>30    \>16\\
                \>                \>31    \>15b\\
                \>                \>32    \>17\\
\end{tabbing}
\vskip - 18 pt
\noindent Table 3. Locations in [10] and [11] of diagrams of species with axes of symmetry.
\end{table}

\smallskip

\noindent {\bf 5. Doubling}

\noindent A natural question is which designs can be doubled as the design of 10-85-1 of Figure 9 and species $36_s$ was doubled to produce the design of Figure 13a and species $35_4$; each strand is replaced by a pair of strands with the same behaviour.
Doubling is something that weavers actually do, but it was introduced into the weaving literature by Gr\"unbaum and Shephard in their [2], p.~154, applied to square satins like 5-1-1.
Box weave 4-3-1 is plain weave doubled.
It is easy to see that, among the fabrics under consideration here, species are converted into species: $34 \rightarrow 33_4$, $36_1 \rightarrow 35_3$, $36_2$ and $36_s \rightarrow 35_4$, and $39 \rightarrow 38$.
Fabrics of species that are the images here, of level 2 above the originating design, as well as those of levels 3 and 4 that are not such images, namely $33_3$ and 37, cannot themselves be doubled to form isonemal fabrics because their lattice units would be of levels 5 and 6, which are too high.
\begin{Thm}{Isonemal fabrics of species $34$, $36$, and $39$ can be doubled and remain isonemal as fabrics of species $33$, $35$, and $38$, respectively, and those of species $33$, $35$, $37$, and $38$ cannot.}
\end{Thm}
\noindent {\bf 6. Halving}

\noindent One would expect that what had been doubled could be halved although the reverse might not be true since the operations are not really inverses.
Halving is the removal of every other strand in each direction, making an intermediate construction called a pseudofabric by Gr\"unbaum and Shephard [5], p.~25, and then the widening of all the strands uniformly to produce another prefabric.
So halving undoes doubling, but doubling does not normally reverse halving.
These prefabrics do not behave in this respect like those with symmetry axes [11].
Considerations are more detailed but less work.
The discussion in [10], \S 7 first paragraph, is basic.

Since what has been doubled, if halved, returns to how it began, there are whole classes of designs of species $33_4$, $35_3$, $35_4$, and 38 that can be halved and remain isonemal because their designs are coloured in blocks of four.
To consider whether there are species {\it all} of whose members can be halved and remain isonemal, we must consider what becomes of the prefabric symmetries on halving.
All that matters is which symmetries act on cells assigned a single number in the numbering of the cells of blocks of four like the quadrants of the Cartesian plane as in [10], since they are all that is left on halving.
For instance, a quarter-turn at a cell corner permutes the numbers---and not just the four adjacent numbers---and so disappears from $G_1$ of the halved prefabric.
Likewise cell-corner and cell-side half-turns permute numbers and so disappear.
Only operations at cell centres are preserved, as follows.
For an $i$-cell (that is, a cell numbered $i$), let $j$-cells be the cells assigned the other number of the same parity and $k$-cells and $\ell$-cells be the cells of the opposite parity.

A cell-centre quarter-turn in an $i$-cell

\noindent 1. is a cell-centre quarter-turn for all $i$-cells,

\noindent 2. becomes a cell-corner quarter-turn for $j$-cells,

\noindent 3. squared becomes a cell-side half-turn for $k$-cells and $\ell$-cells.

\noindent That is, such a quarter-turn would have one of those effects if only cells numbered $i$, $j$, $k$, or $\ell$ were to remain after the halving. 
This is useful information, but it needs to be used in this reformulation.
A cell-centre quarter-turn will affect $i$-cells

\noindent 1. as a cell-centre quarter-turn if it is in an $i$-cell,

\noindent 2. as a cell-corner quarter-turn if it is in a $j$-cell,

\noindent 3. squared as a cell-side half-turn if it is in a $k$-cell or $\ell$-cell.

Further, a cell-centre half-turn in an $i$-cell

\noindent 1. is a cell-centre half-turn for all $i$-cells,

\noindent 2. becomes a cell-corner half-turn for $j$-cells,

\noindent 3. becomes a cell-side half-turn for $k$-cells and $\ell$-cells.

\noindent  
This useful information needs to be used in the following formulation.
A cell-centre half-turn will affect $i$-cells

\noindent 1. as a cell-centre half-turn if it is in an $i$-cell,

\noindent 2. as a cell-corner half-turn if it is in a $j$-cell,

\noindent 3. as a cell-side half-turn if it is in a $k$-cell or $\ell$-cell.

Because centres of rotation at cell corners disappear, we ignore prefabrics with symmetry groups at levels 3 and 4 and of species 34 and $36_2$ at level 2, which have nothing but cell-corner centres of rotation.
We focus attention on what is left, lattice units with cell-centre corners, of species $36_1$ or 39 at level 1 and species $36_s$ at level 2, beginning at level 1.
Consider a species-$36_1$ or species-$39$ fabric with level-1 lattice unit with a corner in an $i$-cell.
The horizontal and vertical sides of the right triangle of which the lattice unit side is hypotenuse are odd and even, $M_1$ and $N_1$, causing the other end of the side to be in a $k$-cell or an $\ell$-cell.
The end of the side of the next lattice unit in the same oblique direction, however, is again in an $i$-cell.
These two sides together are the side of a lattice unit at level 3, whose centre is similarly in a $j$-cell.

The level-3 lattice unit that arises in this way has quarter-turns in $i$-cells at its corners, in a $j$-cell at its centre, and at $k$-cells and $\ell$-cells at mid-sides.
On $i$-cells these act respectively as cell-centred quarter-turns, as cell-corner quarter-turns, and, when squared, as cell-side half-turns \dia when the cells numbered $j$, $k$, and $\ell$ have been discarded to halve the prefabric.
In the halved prefabric they compose a lattice unit at level 1 of $G_1$ of Roth type $36_1$ (Roth type 39 needs \diabb ).

Consider a level-2 lattice unit with a corner in an $i$-cell.
The sides of the right triangle of which the lattice unit side is hypotenuse are both odd, $M_2$ and $N_2$, causing the other end of each adjacent side to be in a $j$-cell.
Passing on around the unit's boundary, the opposite corner is again in an $i$-cell.
The diagonals of this unit are the sides of different lattice units at level 3 related like chain mail, one having $i$-cells at its corners and $j$-cell at its centre and the other having $j$-cells at its corners and a $i$-cell at its centre.
At mid-sides both have centres of quarter-turns in the central $k$-cell or $\ell$-cell of the level-2 lattice unit.
The previous paragraph now applies to the level-3 lattice unit with the corners in $i$-cells.
The result has type-$36_1$ symmetry.

In the previous paragraph, the reasoning of its previous paragraph would have produced a lattice unit at level 4 corresponding to a perfectly legitimate group of symmetries of the halved prefabric at level 2, but that would be a subgroup of the symmetry group determined by the level-3 lattice unit.

There is a hazard in the apparent conclusion to the above argument, namely that the symmetry group of the halved prefabric has only been shown to {\it contain} the group at level 1 of type $36_1$.
It can be larger; at the extreme the halved prefabric can be the trivial prefabric 1-0-1*, which results from two of the four halvings of 10-85-1 of Figure 9.
The other two halvings are the $(5, 3)$ satin.
What has been shown by the above therefore must be less than one might have hoped.
\begin{Thm}{Isonemal fabrics of species $36_1$, $36_s$, and $39$ can all be halved in all four possible ways to produce isonemal prefabrics whose symmetry groups contain a group of Roth type $36_1$.}
\end{Thm}
Since prefabrics of species $36_1$ are carried to species $36_1$, one might wonder whether any fabric has itself as an image.
The odd square satins map to themselves.

\smallskip
\noindent {\bf 7. Fabrics of a given order}

\smallskip
\noindent As in [10] and [11], information in sections 3 and 4 allows all fabrics of a given order to be found because the symmetry groups for a given order can be determined.
The relations among the various numbers in section 3 above are summarized in Table 1.
If fabric designs with quarter-turn symmetries of a given order greater than 4 are wanted, then there may be several different species (or none) where they can be found.
There are two jointly sufficient conditions for existence.
The potential order can be factored into a power $2^p$ and an odd factor $f$.
The first condition is that $f$ is a sum of relatively prime squares, a condition that can be satisfied non-uniquely; $64 + 1 = 49 + 16$.
The second condition is that the power $2^p = 1, 2$, or 4 from column 5 of Table 1.
If these conditions are satisfied, then each relatively prime sum-of-squares decomposition of $f$ can be used to form a level-1 lattice unit of area $f$ that can be used for $G_1$ for designs of species $36_1$ if $p = 0$.
If $p = 1$, then the level-1 unit can be used for designs of species 39, the corresponding lattice unit at level 2 can be used for designs of species 34, $36_2$, and $36_s$, and the corresponding lattice unit at level 3 can be used for designs of species $33_3$ and $35_3$.
If $p = 2$, then the corresponding lattice unit at level 3 can be used for designs of species 38, and the corresponding lattice unit at level 4 can be used for designs of species $33_4$, $35_4$, and 37.
This last species 37, with designs of pure genus V, is the only one in which the phenomenon of falling apart can occur.
Prefabrics in the other species with rotational symmetry are fabrics, as Roth shows in his [8], although not in these terms.
So when the species is 37, prefabrics that fall apart can either be created in the course of working through all and set aside or be specially created and then avoided in the production of fabrics.
For each useful $G_1$, the number of independent cell orbits can be determined and then all of the possibilities examined, discarding prefabrics of species 37 that fall apart and designs with groups of which the desired group is a subgroup as in [10] and [11].
\begin{Thm}{A family of prefabric designs with rotational symmetry and order $f2^p$ greater than four with $f$ odd exists if and only if $f$ is an odd sum of relatively prime squares and $p$ is 0, 1, or 2.}
\end{Thm}
The results of sections 5--7 complete, with the results of [10] and [11], answers for all isonemal prefabrics to their three questions.

\smallskip
\noindent {\bf 8. Woven cubes}

\noindent A topic in the literature to which the taxonomy of Roth is applicable is the cubic case of Jean Pedersen's [6] general problem of woven polyhedra (cf. [1]).
The boundary nuisance that is ordinarily avoided by working in an infinite plane $E$ can be avoided instead by working in a net for a polyhedron, in particular here the six faces of a cube (Figure 16a).
This is the way in which Pedersen represents cubes ([6], Figures 8, 10, 17), and all that it requires is to focus on a portion of the standard weaving diagram as in her Figure 17 (Figure 16a here), due to G.C.~Shephard, which represents the only one of Pedersen's examples with an oblique lattice unit.
The original illustrates weaving with the much more symmetrical plain weave.
\begin{figure}
\begin{tabbing}
\noindent\epsffile{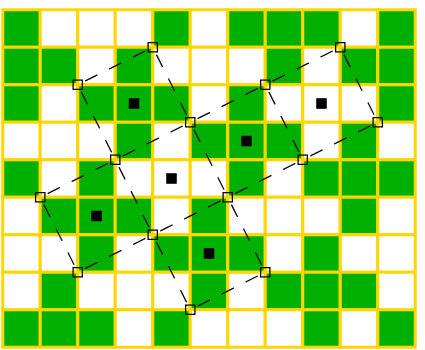}\hskip 10 pt\=\epsffile{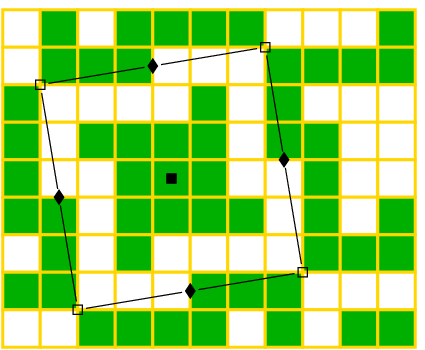}\\

\noindent (a)\>(b)\\
\end{tabbing}
\vskip -18 pt
\noindent Figure 16. a. Shephard's net illustrated in [6], Figure 17, shown here for 10-93-1 with the boundaries of the faces dashed.\hskip 10 pt b. A larger lattice unit for another species-39 fabric.
\end{figure}

As in Figure 16a, we ignore Pedersen's main cubical focus, weaving with lattice units conformed to the directions of the strands (not because they lack interest but because our analysis does not apply to the small set of them).
We also consider only weaving with strands that are---as throughout---perpendicular and such that, at each point of $E$---and the cube net in particular---not on the boundary of a strand, one strand passes over one strand (2-way %
%
%
2-fold prefabrics in the terms of [5]).
We cannot expect that the reflective symmetries of the cube will be symmetries of the strands (consider opposite faces in Figure 16a), and so we leave them out of account.
The six faces of the cube must be woven the same way; in particular, the weaving must be what passes for invariant%
\footnote{Only what passes for invariant because $\tau$ may be involved, as in Figure 16, in the symmetries.}
 under rotation of the cube around its 4-fold, 3-fold, and 2-fold axes, which must collectively be transitive on the strands (Pedersen's definition of isonemality for cubes).
The general demand requires that a face of the net must be related to the other faces by translations or rotations that are symmetries.
The four-fold requirement is that a centre of quarter-turn symmetry must appear in the centre of each face of the net.
The two-fold requirement is that a centre of half-turn symmetry must appear at the mid-sides of each face.
And the three-fold requirement includes that the corners of the faces fall at cell corners; a cell or cell side containing a corner of the net is plainly impossible.
If we require conversely that centres of rotational symmetries of the prefabrics be confined to what will be face centres and corners of the cube and the centres of edges, which we are not bound to do (Pedersen did not), then the net of the cube must be made up of lattice units of the prefabric.
Here we shall tentatively take this further restriction as added to the definition of an isonemal weaving of a cube.
Without the final optional requirement, the match is lost between the symmetries of cube and design.
With these natural restrictions, weaving cubes becomes a matter of choosing designs having the appropriate symmetry groups and therefore determining whole species of designs.

Which lattice units of isonemal prefabrics can be put together in this way?
The sorts of lattice unit that are available are displayed in Figure 10.
That lattice-unit corners fall at cell corners limits either the level to 1 and the centre operation to \blbox when the centre of the unit is at a cell centre as in Figures 10a and 16 (species $36_1$, 39) or the level to 2, 3, or 4 when the centre of the unit is at cell corners as in Figure 10b (species 33, 34, 35, $36_2$, 37, 38).
Elimination of level 2 with lattice-unit centre at cell centre rules out species $36_s$.

Figures 6a and 7 display species of design with \blbox at the centre of level-1 lattice units (and in each case alternative descriptions of the same symmetry groups).
When one thinks about translates of the lower lattice unit in Figure 7 making a net like that in Figure 16a, it looks as though there is a difficulty when the \whbox at the centre of the alternative lattice unit is taken as a corner of the cube requiring 3-fold symmetry: there will be either two dark crosses and one pale cross used for the cube or vice versa.
But this apparent difficulty is a result of the colouring convention; there is no difficulty with the topology, since dark and pale have to do with directions in the diagram and are not inherent in the topology.
A dark cross is the image under a quarter-turn of a pale cross: three rotations of an initial pale cross bring it back to itself as it should because the offending second pale cross is missing from the net.
The colouring convention {\it conceals} the real difficulty in using a net based on the upper lattice unit of Figure 6a.
Because \blbox includes $\tau$, when it is applied three times the face that is the lattice unit rotated is returned to its original position with colours complemented because $\tau$ has been applied three times.
Designs of species 39 are not ruled out by these considerations, but those of species $36_1$ are ruled out.
Species 39 must be used with the lattice unit that is the lower in Figure 7, not the higher.
Since the argument against allowing \blbox at the corner of a lattice unit is general, it gives a lemma.

\begin{Lem}{Lattice units of $G_1$ of an isonemal prefabric cannot compose the net of a woven cube if they have \blbox at their corners.}
\end{Lem}
\noindent The lemma rules out species $35_3$, $35_4$, and $36_2$, leaving, $33_3$ and $33_4$, which differ only in level, 37 and 38, which likewise differ only in level and have the symmetry combination of species 39 but different centre location, 34, and 39.
Glancing at Figures 12a and 14 suggests that the lemma rules out species 37 and 38 too, but one must recall that there are alternative lattice units for these species as for species 39 having \blbox at their centres and the required \whbox at each corner.

The restriction of lattice units to levels 1 to 4 and on the dimensions of them at level 1 were put in place to ensure that the strands in $E$ not fall into classes not related to one another by symmetries of their design.
Since the centres of half-turns (including $\hbox{\whboxx\hskip 0.7 pt}^2 = \hbox{\blboxx\hskip 0.7 pt}^2$) relating parallel strands can be far apart in $E$ and still perform their task, farther apart than would keep them inside a cubical net, one wonders whether all the lattice units of species 33, 34, 37, 38, and 39 can compose cubical nets in which the rotational isometries of the cube are transitive on the strands.

\begin{Thm}{A necessary and sufficient condition that lattice units of an isonemal prefabric of order greater than 4 be the regions of a net of an isonemal woven cube is that the prefabric be of species $33$, $34$, $37$, $38$, or $39$.}
\end{Thm}

\noindent {\it Proof.} The definition of isonemal woven cube rules out species with more symmetry than those listed, those like plain weave whose symmetry groups contain groups of the listed Roth types as proper subgroups.

Necessity. If the design of one face of an isonemal woven cube is set out on the plane and propagated across the plane in accordance with its own symmetry group, including centres of half-turns at mid-sides and of side-preserving quarter-turns at its corners, it will cover the plane to be sure.
Its images will be the lattice units of an isonemal prefabric with one of the listed groups.
At the face centre will be a quarter-turn centre, either \blbox or \whboxx .
If it is a lattice unit at level 1, then the Roth type of the group is 39.
If it is at level 2, then the Roth type is 34.
If it is at level 3, then the Roth type is $33_3$ or 38.
If it is at level 4, then the Roth type is $33_4$ or 37.
It cannot be of a higher level, since the isonemality of the cube demands that there be strand-to-adjacent-strand transformations, which must be half-turns,  quarter-turns, or translations, and at a higher level there are not sufficient centres even in the whole of $E$ for half-turns or quarter-turns to relate all adjacent strands, nor at higher levels can translations do it.

Sufficiency. Lattice units of an isonemal prefabric design of the listed species obviously fit together to form nets of cubes, and there are transformations from strands to perpendicular strands. 
What needs to be proved is that, if the design is of one of the listed species, there are transformations from strands to adjacent parallel strands.
This is obviously the case if centres of rotation are allowed to pass outside the net by the isonemality of the prefabric, that is, one strand is mapped to each adjacent strand by a half-turn, which may be the square of a quarter-turn.
We must allow for only what happens on the surface of the cube.
Consider two strands adjacent over some distance; we shall see that there is at least one centre of rotation on their common boundary and within a net-like covering of the cube and so with image on the cube.

First consider such a pair with a common boundary $B$ along their whole length, which being in the cube is finite.
In $E$, their common boundary is an infinite straight line $L$ along which there are periodically centres of rotation interchanging them by half-turns.
When one thinks of a portion of $E$ large enough to contain, say, two such centres of rotation as a net-like multiple covering of the cube, one sees that $L$ is mapped round and round the finite $B$.
Accordingly, the centres of rotation on $L$ must be mapped to at least one centre of rotation along $B$.

\begin{figure}\noindent\epsffile{19.eps}

\noindent Figure 17. Net-like diagram of a cube weavable, for instance, with the fabric of Figure 16b showing the common boundary of two strands and the centre of half-turn at its mid-point.
\end{figure}

Secondly, consider two strands adjacent over some distance that do not behave so conveniently.
They fail to do so if their common boundary $B$ ends in each direction at a vertex of the cube (then the two strands briefly go their separate ways); otherwise they behave as above.
If, starting at some point on $B$ one passes along it and hits a vertex then one must hit a vertex if one passes along $B$ in the opposite direction.
$B$ is the image in the cube of a line segment $\ell$ in $E$ between vertices of lattice units.
Since the lattice units are square, the mid-point of $\ell$ falls on the mid-point of a lattice-unit side or centre point of a lattice unit; if it were to fall on a vertex---the only other possibility, then that vertex would map to an end of $B$, contradicting the definition of $\ell$.
There are centres of half-turns at all of the possible mid-points, and so there is a half-turn relating the two strands at the mid-point of $B$.
As this is the case for all adjacent parallel strands, the cube is woven isonemally.

The theorem is proved.

Examples will make clearer what happens in the more interesting of the two cases discussed in the proof of sufficiency.
The kind of line segment $\ell$ that maps on $B$ in Figure 16a appears there three times horizontally beginning at the top left, top right, and bottom right corners of the lattice unit of the central pale cross and two times vertically from the top right and bottom left corners of its lattice unit. 
The mid-point of each of these five line segments is the middle of an edge of the cube, where there is a half-turn centre.

Any similar covering of a cube can support the same reasoning.
Figure 17 is of a covering of a cube where the top and bottom are not shown opposite each other to save space and there are six squares to map to the faces of the cube, two faces being covered twice.
This is so as to accommodate the whole $\ell$ mapping to the boundary $B$, which is the horizontal line in Figure 17.
It is like those in Figure 16a, but it goes around the cube 1 1/2 times instead of half-way around.
If one takes the lattice unit of Figure 16b as a face of the cube and $\ell$ as beginning in the lattice unit's top left corner, then it returns to that lattice unit two-thirds of the way down its left side and goes out the right side again, needing another face before it reaches its other vertex endpoint.
As Figure 16b indicates, such a boundary, being a strand edge, sinks by one strand each time it crosses a face.
This can be seen also in Figure 17. 
The boundary must therefore pass across six faces to fall from the top of its first face to the bottom of its sixth.
As Figure 17 shows, $\ell$ has a half-turn centre at its mid-point.
This it encounters on entering the fourth of the six faces it crosses, the one in which it began being the fifth.

The strands above and below $\ell$ through the middle of the diagram have it as their common boundary over its length but then go their separate ways beyond the single cell beyond each end of it, which they pass through and beyond in perpendicular directions.
It is an interesting exercise to trace the two paths; each runs through seven lattice units to get from the right side of the diagram to the left, the end ones being common to both strands.

\medskip
\noindent {\bf References}

\smallskip

\noindent [1] Paulus Gerdes: Weaving polyhedra in African cultures. {\it Symmetry: Culture and Science} {\bf 13} (2002), 339--355.

\noindent [2] Branko Gr\"unbaum and Geoffrey C.~Shephard: Satins and twills: An introduction to the geometry of fabrics, {\it Mathematics Magazine} {\bf 53} (1980), 139--161.

\noindent [3] ------: A catalogue of isonemal fabrics, in {\it Discrete Geometry and Convexity}, Jacob E.~Goodman {\it et al.,} eds.~{\it Annals of the New York Academy of Sciences} {\bf 440} (1985), 279--298.

\noindent [4] ------: An extension to the catalogue of isonemal fabrics, {\it Discrete Mathematics} {\bf 60} (1986), 155--192.

\noindent [5] ------: Isonemal fabrics, {\it American Mathematical Monthly} {\bf 95}(1988), 5--30.


\noindent [5] J.A.~Hoskins and R.S.D.~Thomas: The patterns of the isonemal two-colour two-way two-fold fabrics, {\it Bull. Austral.~Math.~Soc.}~{\bf 44} (1991), 33--43.

\noindent [6] Jean Pedersen: `Some isonemal fabrics on polyhedral surfaces', in {\it The geometric vein}, pp.~99--122. Chandler Davis, Branko Gr\"unbaum, and F.A.~Sherk, eds.~New York: Springer, 1981.

\noindent [7] Richard L.~Roth: The symmetry groups of periodic isonemal fabrics, {\it Geometriae Dedicata} {\bf 48} (1993), 191--210.

\noindent [8] ------: Perfect colorings of isonemal fabrics using two colours, {\it Geometriae Dedicata} {\bf 56} (1995), 307--326.

\noindent [9] D.~Schattschneider: The plane symmetry groups: Their recognition and notation, {\it American Mathematical Monthly} {\bf 85} (1978), 439--450.

\noindent [10] R.S.D.~Thomas: Isonemal prefabrics with only parallel axes of symmetry, {\it Discrete Mathematics} {\bf 309} (2009), 2696--2711. 
Available at arxiv.org/abs/math/0612808v2.

\noindent [11] ------: Isonemal prefabrics with perpendicular axes of symmetry, {\it Utilitas Mathematica}, to appear. Available at arxiv.org/abs/0805.3791v1.

\end{document}